 \newif \ifSIAM
 \documentclass[10pt,a4paper]{article}
 \usepackage{authblk}
 \usepackage{graphicx}
 \usepackage{amssymb}
 \usepackage{amsmath}
 \usepackage{amsthm}
 \usepackage{booktabs}
 \usepackage{paralist}
 \usepackage{bm}
 \usepackage{url}
 \usepackage{color}
 \usepackage{mathtools}
 \usepackage{subcaption}
  \usepackage{array}
  \usepackage{multirow}
  
 \newcolumntype{M}[1]{>{\centering\arraybackslash}m{#1}}
 
 \numberwithin{equation}{section}

 \newif \ifbbm
 \newif \ifPDF
 \PDFtrue

\newtheorem{lemma}{Lemma}[section]

\newtheorem{proposition}[lemma]{Proposition}

 \newcommand\Om{\Omega}
 \newcommand\crl{\mathrm{curl}}
 \newcommand\dv{\mathrm{div}}
 \newcommand\Hc{{H}(\crl, \Om)}
 \newcommand\Hcz{{H}_0(\crl, \Om)}
 \newcommand\R{\mathcal{R}}
 \newcommand\LL{\mathcal{L}}
 \newcommand\RR{{\mathbb R}}
 \newcommand\tu{\bm{u}}
 \newcommand\tv{\bm{v}}
 \newcommand\tn{\bm{\mathrm{n}}}
 \newcommand\tE{\bm{E}}
 \newcommand\tH{\bm{H}}
 \newcommand\tJ{\bm{J}}
 \newcommand\bphi{\bm{\phi}}
 \newcommand\bpsi{\bm{\psi}}
 \newcommand\tx{{\bm{x}}}

 \newcommand\F{\mathcal{F}}
 \newcommand\I{\mathcal{I}}
 \newcommand\te{\bm{e}}
 \newcommand\ii{\mathrm{i}} % index i 
 \newcommand\jj{\mathrm{j}} % index j
 \newcommand\ZZ{{\mathbb Z}}
 \newcommand\tb{\bm{b}}
 \newcommand\tW{\bm{W}}
 \newcommand\tA{\bm{A}}
 \newcommand\N{\mathcal{N}}
 \newcommand{\bal}{{\bm{\alpha}}}

 \title{Deep Fourier Residual method for solving time-harmonic Maxwell's equations}

 \author[1]{Jamie M. Taylor} 
 \author[2]{Manuela Bastidas}
 \author[2,3,4]{David Pardo}
 \author[5]{Ignacio Muga}
 
 \affil[1]{CUNEF University, Madrid, Spain. \protect\\ \texttt{jamie.taylor@cunef.edu}}
 \affil[2]{University of the Basque Country (UPV/EHU), Leioa, Spain.}
 \affil[3]{Basque Centre for Applied Mathematics (BCAM), Bilbao, Spain.}
 \affil[4]{Ikerbasque (Basque Foundation For Sciences), Bilbao, Spain.}
 \affil[5]{Pontificia Universidad Católica de Valparaíso, Chile.}
 
 \ifPDF \usepackage[breaklinks,bookmarks = false]{hyperref}

 \hypersetup{colorlinks, linkcolor = blue, citecolor = blue,
urlcolor = blue, pagecolor = blue, plainpages = false, pdfwindowui = false,
pdfstartview = {FitH}, pdftitle = {DFRHcurl},
pdfauthor = {ManuelaBastidas} }
 \fi

\begin{document}
 \maketitle

 %\tableofcontents
 
 \begin{abstract}
 	Solving PDEs with machine learning techniques has become a popular alternative to conventional methods.  In this context, Neural networks (NNs) are among the most commonly used machine learning tools, and in those models, the choice of an appropriate loss function is critical. In general, the main goal is to guarantee that minimizing the loss during training translates to minimizing the error in the solution at the same rate. In this work, we focus on the time-harmonic Maxwell's equations, whose weak formulation takes $\Hcz$ as the space of test functions. We propose a NN in which the loss function is a computable approximation of the dual norm of the weak-form PDE residual. To that end, we employ the Helmholtz decomposition of the space $\Hcz$ and construct an orthonormal basis for this space in two and three spatial dimensions.  Here, we use the Discrete Sine/Cosine Transform to accurately and efficiently compute the discrete version of our proposed loss function. Moreover, in the numerical examples we show a high correlation between the proposed loss function and the $H$(curl)-norm of the error, even in problems with low-regularity solutions.
 \end{abstract}

 \section{Introduction}
 
 Partial Differential Equations (PDEs) are essential tools for modeling and simulating various scientific and industrial problems, in particular, they form the backbone of modern physics. Herein, we concentrate our attention on Maxwell's equations. This set of equations describes the propagation of electromagnetic waves through different media. Commonly used methods to calculate the solution of Maxwell's problems range from exact methods \cite{gralak2012exact,jones1979methods} to numerical approximations, such as Finite Differences, Finite Elements (FEM), and Discontinuous Galerkin \cite{camargo2020hdg, chaumont2022frequency, kirsch2015mathematical, monk2003finite, taflove2005computational}. Nevertheless, exact solutions can only be obtained in rare scenarios, and proposing a proper numerical technique that is conformal, accurate, and efficient is often challenging. 
 
Using machine learning techniques to approximate the solution of Maxwell's equations is an attractive alternative to classical approaches. Some of the most used machine learning tools in the context of PDEs are Neural Networks (NNs) \cite{calabro2022effect, doleglo2022deep,lecun2015deep,lu2021deepxde,raissi2019physics,ruthotto2020deep}. These architectures have shown promising results when tackling complex nonlinear systems of equations that underlay physical phenomena. For instance, we highlight their potential to solve parametric PDEs \cite{han2018solving,khoo2021solving}, enhance classical numerical methods \cite{brevis2023learning} and solve problems in presence of singularities or sharp gradients \cite{calabro2021extreme, sluzalec2023quasi}. 
 
 Some popular methods like PINNs (Physics-Informed Neural Networks) \cite{karniadakis2021physics, raissi2019physics} enforce a NN to satisfy the strong formulation of a PDE by implementing a numerically tractable norm of the strong-form residual as a loss function. Even though PINNs can approximate the solutions of many physical problems, they are inaccurate when the weak solution does not satisfy the strong-form equation. As a result, PINNs are severely constrained in many applications that naturally produce low-regularity solutions. Maxwell's equations, for example, may admit a solution in the functional space $\Hc\setminus H^2(\Om)$ (or even in $\Hc\setminus H^1(\Om)$ as shown in, e.g., \cite{rukavishnikov2012new}) since the smoothness of the solution depends on the regularity of the domain, the sources, and the boundary conditions. In these situations, PINNs may not be applicable.
 
 In contrast, VPINNs (Variational Physics-Informed Neural Networks) \cite{kharazmi2019vpinns, kharazmi2021hp} use the residual of the weak formulation of a PDE in the loss function. In this approach, it is essential to select an appropriate set of test functions and define a computable loss that controls the error of the solution, which is generally non-trivial.  
 
An ideal loss function would be the energy norm of the error. Since the error function is unavailable in practice (one would need the exact solution), one typically resorts to minimizing the dual norm of the weak-form PDE residual, that is $\R: H \to H'$, where $H$ is a Hilbert space of test functions and $H'$ its dual. 
 In many practical cases of linear PDEs, if $u^*$ is the exact solution of a PDE, then, there exist constants $0 < \gamma < M$, such that 
 \begin{equation}\label{eq:gbound}
 	\frac{1}{M} \| \R(u) \|_{H'} \leq \|u- u^*\|_H \leq \frac{1}{\gamma} \| \R(u) \|_{H'},
 \end{equation}
 where $u$ is an approximation to the solution $u^*$ \cite{taylor2023deep}. It is clear from \eqref{eq:gbound} that the dual norm of the residual $\| \R(u) \|_{H'}$ is equivalent to the energy norm of the error, and so it is an appropriate choice of loss function, i.e., $\LL(u)  := \| \R(u) \|_{H'}$. Unfortunately, it is in general very challenging to evaluate $\| \cdot \|_{H'}$. 
 
 For similar problems based on $H^1$ test function spaces, the authors of \cite{taylor2023deep} proposed the Deep Fourier Residual (DFR) method. They employed a numerical method for approximating the dual norm of residuals corresponding to PDEs with $H^1$ test function spaces via a spectral representation of the dual norm, which may be implemented using the Fast Fourier Transform (FFT). 
 
 This paper follows the same spirit; we extend the ideas of \cite{taylor2023deep} and numerically implement the dual norm as a loss function to solve PDEs with $\Hcz$ test function spaces, motivated by the time-harmonic Maxwell's equations. Similarly to the DFR method for PDEs with $H^1$ test function spaces, here the key challenge is to construct an appropriate orthonormal basis for $\Hcz$. To find such a basis, we use the classical Helmholtz decomposition of a vector field \cite{monk2003finite} and construct a complete set of basis functions for the space $\Hcz$ on product domains in two and three spatial dimensions. The strength of the DFR method lies in solving problems with low regularity, where PINNs-like methods based on strong formulations fail. Moreover, our choice of loss function bounds the $H$(curl)-norm of the error of the solution. In this way, minimizing the loss during training implies the reduction of the error solution at the same rate. We provide numerical examples that demonstrate the correlation between the proposed loss function and the $H$(curl)-norm of the error.  In the specific case of Maxwell's equations, the DFR method produces accurate results on heterogeneous and discontinuous media, as well as strong correlations between the loss and $H$(curl)-norm of the error during training. 
 
Notice that, besides our work, other studies employ the combination of Fourier basis functions and machine learning techniques. This is usually called Physics-informed Spectral Learning (PiSL) \cite{espath2023physics,espath2021statistical}. Nevertheless, these techniques are associated with data analysis and in contrast to our study, they are situated within the framework of physics-informed statistical learning.
 
 Our proposed DFR method faces significant challenges that are similar to those encountered in \cite{taylor2023deep}. Here, we only consider product domains, to profit from the results in \cite{costabel2019maxwell}, and with Dirichlet-type boundary conditions. Different strategies must be used to tackle problems with non-trivial geometries and general boundary conditions. Specifically, using the DFR method when the problem involves general geometries is challenging since constructing an orthonormal basis for the space $\Hcz$ is not straightforward. Another important limitation of our technique is that the constants bounding the norm of the error, in \eqref{eq:gbound}, diverge in certain cases, e.g., when one considers materials with a large variation in parameters or frequencies close to a resonant frequency of the system. These technical complications are seen as stability issues intrinsic to the PDE similar to those encountered in traditional approaches such as FEM \cite{demkowicz1998modeling}, where the accuracy of the error estimates deteriorate under the same conditions.
 
 The remainder of this work is organized as follows. In Section~\ref{sec:problem}, we state Maxwell's equations and the weak formulation of the problem. There, we also motivate the choice of the dual norm of the residual operator as a natural loss function for training NNs in the case of Maxwell's equations. Later, we discuss the relevance of constructing an appropriate orthonormal basis for the space $\Hcz$. The details of generating such an orthonormal basis, using Helmholtz decomposition in two and three spatial dimensions, are described in Section~\ref{sec:basisF2}, with the technical details deferred to Appendix~\ref{sec:appA}. In Section~\ref{sec:basisF2}, we also present the test functions in 2D and 3D for the simple case of squared and cubic domains and the corresponding calculations are presented in the Appendix~\ref{sec:appB}. In Section~\ref{sec:dfrNN}, we describe the architecture of a NN and explain the core of the DFR method, which consists of constructing a computable discretized loss function. Finally, in Sections~\ref{sec:numerical} and \ref{sec:concl}, we present numerical experiments, conclusions, and directions for future research.
 
 \section{Problem statement}\label{sec:problem}

Consider a domain $ \Om \subset \RR^n$, with $n = 2$ or $3$, whose boundary $ \Gamma  := \partial \Om$ is polyhedral and connected. Given an impressed field $ \tE^{I}$, and an electric density current source $ \tJ$, we look for an electric field $ \tE$, and a magnetic field $ \tH$ solving the so-called macroscopic linear Maxwell's equations in a time-harmonic form
\begin{equation}\label{eq:maxwell}
	\begin{aligned}
		\crl(\tE) - i \omega \mu \tH 		 & = 0 	 & &\text{in } \Om,	 \qquad \text{(Faraday's Law)}, \\
		\crl (\tH) + i \omega \epsilon \tE 	& = \tJ & &\text{in } \Om, \qquad \text{(Ampere's Law)}, \\
		\tE \times \tn 	& = \tE^{I} & &\text{on } \Gamma,
	\end{aligned} 
\end{equation}
where $i$ is the imaginary unit, $\tn$ denotes the outward unit normal vector of $\Om$, $\omega \in \RR$ is the angular frequency, and $\mu$ and $\epsilon$ are space-dependent functions standing for the magnetic permeability and electrical permittivity, respectively. These functions generally may be tensor-valued, but we consider here only the scalar case.

The formulation  \eqref{eq:maxwell} arises by considering the time-dependent Maxwell's equations under the following Ansatz on the electric and magnetic fields
\begin{equation*}
	\tE(x,t) = \mathfrak{R}(e^{i \omega t}\tE(x)), \text{ and } \tH(x,t) = \mathfrak{R}(e^{i \omega t}\tH(x)).
\end{equation*}

In \eqref{eq:maxwell}, the dependency of the time $t>0$ is implicit. The model \eqref{eq:maxwell} is completed with the following Gauss' Laws for magnetic and electric fields
\begin{equation}\label{eq:diveq}
	\begin{aligned}
		\dv (\mu \tH) 		& = 0 		& \text{in } \Om, \\
		\dv (\epsilon \tE) & = \rho & \text{in } \Om,
	\end{aligned} 
\end{equation}
where $\rho$ is the density of free charge. 

For the sake of simplicity, in this work we only consider the case when $\omega$ is non-zero. Notice that when $\omega \neq 0$ the equations in \eqref{eq:diveq} are consequences of \eqref{eq:maxwell} and the continuity equation 
\begin{equation*}
	i \omega \rho + \dv (\tJ) = 0, 
\end{equation*}
which relates the rate of change of the charge density to the divergence of the current density.

We point out that we have only considered Dirichlet-type boundary conditions. Other strategies must be implemented for different boundary conditions. Whilst the following discussion will only be limited to the case of homogeneous boundary conditions, as our approach is based on an analysis on the space of test functions, rather than trial functions, the extension to non-homogeneous boundary conditions is straightforward as the test function space remains unchanged.

\subsection{Preliminaries}
We start by defining relevant operators and function spaces. Herein, the definitions are regarded as classical and can be found in \cite{girault2012finite}. First, we denote $\mathcal{D}(\Om)$ the space of smooth functions with compact support in $\Om$ and $\mathcal{D}'(\Om)$ is the space of distributions. We let $L^p(\Om)$ be the space of the $p$-integrable real-valued functions equipped with the usual norm. 

For smooth functions, we define the divergence, gradient and curl according to their usual definitions (see \cite{schey2005div}). In 2D, we define the curl and its adjoint as 
\begin{equation}\label{curl2}
	\begin{aligned}
		\crl : 		& \quad [\mathcal{D}(\Om)]^2 \ni \bphi \quad \mapsto \quad \partial_y\bphi_1-\partial_x\bphi_2 \in\mathcal{D}(\Om),\\
		\crl^*:	 & \quad \mathcal{D}(\Om) \ni \phi \quad \mapsto \quad \begin{pmatrix}-\partial_y\phi\\ \partial_x\phi \end{pmatrix} \in [\mathcal{D}(\Om)]^2.
	\end{aligned}
\end{equation}
The curl operator in 2D and 3D can then be extended as mappings from appropriate $L^2$ spaces to $\mathcal{D}'(\Om)$ by duality. We then define the function space $\Hc$ consisting of functions in $[L^2(\Omega)]^n$ whose curl, interpreted in the sense of distributions, is in $[L^2(\Omega)]^{n'}$, i.e.,
\begin{equation*}
		\Hc  :=  \{\tu\in [L^2(\Om)]^n : \crl(\tu)\in [L^2(\Om)]^{n'}\},
\end{equation*}
where $n' = 1$ if $n = 2$ and $n' = 3$ if $n = 3$. The space $\Hc$ is a Hilbert space with inner product given by
\begin{equation}\label{eq:innerprod}
		( \tu, \tv )_{\Hc}  = \int_\Om \crl(\tu)\cdot\crl(\tv) + \tu \cdot \tv \, d\tx \qquad \forall \tu, \tv \in \Hc.
\end{equation}

For any bounded Lipschitz domain $\Om \subset \RR^n$ with boundary $\Gamma$ and outward normal $\tn$, the mapping $\gamma_t: \mathcal{C}^1(\bar{\Om}) \to L^2(\Gamma)$ with $\gamma_t(\tu) = \tu|_{ \Gamma} \times \tn$ can be uniquely extended to the continuous tangential trace operator, $\gamma_t:\Hc\to H^{-\frac{1}{2}}( \Gamma,\RR^{d'})$ (see \cite{buffa2002traces})\footnote[1]{In two dimensions, the cross product is interpreted as the scalar product $ \tu \times \tv = \tv_2 \tu_1 - \tu_2 \tv_1$.}. We then define the space
\begin{equation*}
	\begin{split}
		\Hcz  := &\{ \tu\in \Hc: \, \gamma_t(\tu) = 0\}.
	\end{split}
\end{equation*}

\subsection{Weak formulation}

There are multiple weak formulations for the Maxwell system, all based on the general idea of minimizing the functional that represents the energy of an electromagnetic field. Notice that in \eqref{eq:maxwell} one can eliminate $\tH$ or $\tE$ from each of the equations. Assuming $\mu$ and $\epsilon$ are real-valued, bounded and non-zero functions, $\tJ \in [L^2(\Om)]^n$, and $\tE^{I} = \mathbf{0}$, the weak formulation corresponding to the electric field in the problem \eqref{eq:maxwell} is: Find $\tE \in \Hcz$ satisfying
\begin{equation}\label{eq:weakMax}
	\int_\Om \mu^{-1} \crl(\tE) \cdot \crl(\bphi)- \omega^2 \epsilon \tE \cdot \bphi \, d\tx = \int_\Om i \omega \tJ \cdot \bphi \, d\tx \qquad \forall \bphi \in \Hcz.
\end{equation}

An analogous weak form exists for the magnetic field $\tH$. Notice that, Gauss' Laws in \eqref{eq:diveq} are satisfied weakly, by considering test functions $\bphi = \nabla u$ for $u \in H^1_0(\Om)$ in \eqref{eq:weakMax}. 

\subsection{Residual minimization}

The residual operator corresponding with the weak form \eqref{eq:weakMax} is $\R : H \to H'$ with $H = \Hcz$ and $H'$ being its dual. This weak residual operator may be expressed in the general form 
\begin{equation}\label{eq:maxshort}
	\langle \R(\tE), \bphi \rangle_{H'\times H} = b(\tE,\bphi) - \ell(\bphi),
\end{equation}
where $\ell \in H'$ is 
\begin{equation*}
	\ell(\bphi) = \int_\Om i \omega \tJ \cdot \bphi \, d\tx\qquad \forall \bphi \in H ,
\end{equation*}
and $b: H \times H \to \RR$ is the following bilinear form
\begin{equation}\label{eq:bilin}
	b(\tE,\bphi) = \int_\Om \mu^{-1} \crl(\tE) \cdot \crl(\bphi)- \omega^2 \epsilon \tE \cdot \bphi \, d\tx \qquad \forall \bphi \in H.
\end{equation}

The existence and uniqueness of a solution of \eqref{eq:weakMax}, for $\omega$ outside of a countable set of resonant frequencies, is proved in, e.g. \cite[Chapter 4]{monk2003finite} and \cite[Theorem~4.32]{kirsch2015mathematical}. Moreover, using the reasoning in  \cite[Section~25.3]{ern2021finite}, we know that the solution for the variational problem exists and is unique if and only if the following bounds apply 
\begin{equation}\label{eq:bounds}
	\begin{aligned}
		\| \R(\tE) \|_{H'}^2 & = \sup_{\bphi \in H \backslash \{ 0 \}} \frac{| \langle \R(\tE), \bphi \rangle_{H'\times H} |}{ \| \bphi \|_{H}} &\leq M \| \tE - \tE^* \|_{H} ,\\
		\| \R(\tE) \|_{H'}^2 & = \sup_{\bphi \in H \backslash \{ 0 \}} \frac{| \langle \R(\tE), \bphi \rangle_{H'\times H} |}{ \| \bphi \|_{H}} &\geq \gamma \| \tE - \tE^* \|_{H},
	\end{aligned}
\end{equation}
where $ \gamma$ and $M$ are positive constants depending on $\mu$, $\omega$ and $\epsilon$, and $\tE^*$ denotes the exact solution of \eqref{eq:weakMax}. We emphasize that the coercive case of $\epsilon < 0$, is mathematically interesting, as it implies that the bilinear form \eqref{eq:bilin} becomes equivalent to the inner product on $\Hc$. 

On the other hand, for any $\F \in H'$, by the Riesz representation theorem, there exists some $\tu_{\F} \in H$ with $\F(\tv) = (\tu_{\F} ,\tv)_H$ for all $\tv\in H$ and $ \| {\F} \|_{H'} = \| \tu_{\F} \|_{H}$. Furthermore, if $(\Phi_k)_{k \in \I}$ is an orthonormal basis of $H$, with $\I$ denoting a set of indices, then, by using the generalized Parseval's identity, we have that the dual norm of any $ {\F} \in H'$ can be expressed as
\begin{equation}\label{eq:parse}
	\| {\F} \|_{H'}^2 = \| \tu_{\F} \|_{H}^2 \stackrel{\text{Parseval}}{ = } \sum\limits_{k\in\I} ( \tu_{\F} ,\Phi_k )_H^2 = \sum\limits_{k\in\I} {\F}(\Phi_k)^2.
\end{equation}

From \eqref{eq:maxshort} and \eqref{eq:parse} we obtain an expression for the dual norm of the residual as 
\begin{equation}\label{eq:parseMax}
	\| \R(\tE) \|_{H'}^2 = \sum\limits_{k\in \I} \langle \R(\tE), \Phi_k \rangle_{H'\times H}^2, 
\end{equation}
 According to \eqref{eq:parseMax}, determining a set of orthonormal basis functions $\Phi_k$ is all that is required to calculate the dual norm of the residual. Whilst this is generally a non-trivial task, in the following section, we will construct such a set of basis functions in simplified geometries. 

\section{A set of basis functions for $\Hcz$}\label{sec:basisF2}

In order to find basis functions for $H_0(\crl,\Om)$, we seek an orthonormal eigenbasis for the differential operator corresponding to the inner product \eqref{eq:innerprod}. More details of these ideas can be found in \cite[Chapter 4]{kirsch2015mathematical}. The inner product \eqref{eq:innerprod} is naturally associated with the differential operator $(1+\crl \text{-} \crl)$. So, we consider the weak eigenpairs $(\lambda_k,\Phi_k) \in \RR\times \Hcz$ solving
\begin{equation}\label{eq:eigenCurl0}
  \begin{aligned}
    \int_\Om \crl(\Phi_k)\cdot\crl(\tv) + \Phi_k \cdot \tv \, d\tx & = \lambda_k \int_\Om \Phi_k \cdot \tv \, d\tx \qquad \forall \tv\in \Hcz, \\
    \| \Phi_k \|_{\Hc} & = 1.
  \end{aligned}
\end{equation}
In strong form, we may write 
\begin{equation}\label{eq:eigenCurl}
  \begin{aligned}
    (1+\crl \text{-} \crl) \Phi_k & = \lambda_k \Phi_k &\text{ in } \Om, \\
    \Phi_k \times \tn & = 0 &\text{ on } \Gamma.
  \end{aligned} 
\end{equation}

If the inverse of $(1+\crl \text{-} \crl)$ were to be compact and self-adjoint on $\Hcz$, then, the application of the Hilbert–Schmidt theorem to the inverse would provide the existence of an eigenbasis for the operator itself and a complete eigenbasis of $\Hcz$.
Nevertheless, we notice that, when restricted to functions of the form $\tv = \nabla u$ for $u\in H_0^1(\Om)$, the inverse of the operator $(1+\crl \text{-} \crl)$ is the identity and thus not compact. However, this technical issue can be resolved using the Helmholtz decomposition of the space $\Hcz$, i.e., we decompose the problem \eqref{eq:eigenCurl0} into two sub-problems. We omit the specifics of this reasoning and refer to \cite[Chapter 6]{brezis2011functional} and \cite{davies1995spectral} for more details.

From \eqref{eq:eigenCurl0}, one obtains that the eigenvalue $\lambda_k = 1$ has an infinite-dimensional eigenspace, that is, the null space of the $\crl$ operator. In simply connected domains, we have that $\crl(\tv) = 0$ implies $\tv = \nabla u$ for some $u\in H^1(\Om)$. Moreover, if $u\in H_0^1(\Om)$, then $\nabla u$ is parallel to the unit normal vector $\tn$ on $\partial \Omega$, which means that $\nabla u\in \Hcz$. Thus, we identify a large space of eigenvectors with eigenvalue $\lambda_k = 1$,  $\nabla H_0^1(\Om) \subset \Hcz$ defined as 
\begin{equation*}
  \nabla H^1_0(\Om)  := \{ \nabla u:u\in H^1_0(\Om)\}.
\end{equation*}
We note that equipping $H^1_0(\Om)$ with the inner product $(u,v)_{H^1_0(\Om)} = \int_\Om\nabla u \cdot\nabla v \, d\tx$, the space $\nabla H_0^1(\Om)$, as a subspace of $\Hcz$, is isometric to $H_0^1(\Om)$, i.e., $(\nabla u,\nabla v)_{\Hc} = (u,v)_{H_0^1(\Om)}$. Consequently, the space $\nabla H^1_0(\Om)$ forms a closed subspace of $\Hcz$ (see \cite[Lemma 4.20]{kirsch2015mathematical}), and we can employ the following orthogonal decomposition of $\Hcz$
\begin{equation}\label{eq:decomp}
  \Hcz = X_0(\Om)\oplus \nabla H^1_0(\Om),
\end{equation}
where $X_0(\Om)  := (\nabla H^1_0(\Om))^\perp$. Notice that for a function $\tv \in X_0(\Om)$, one necessarily has that 
\begin{equation*}
  (\tv,\nabla u)_{\Hc} = 0 \qquad \forall u\in H^1_0(\Om),
\end{equation*}
meaning $\dv(\tv) = 0$ weakly. That is to say that vector fields in $\Hcz$ can be decomposed into two parts: a curl-free component and a divergence-free component. Finding an orthonormal basis for $H_0(\crl,\Om)$ reduces to finding an orthonormal basis for each component. We consider each of these in turn. 

\subsection{A set of basis functions for $\nabla H^1_0(\Om)$}

To find a basis for the space $\nabla H^1_0(\Om)$, we use the fact that the differential operator $\nabla$ defines an isometry between $H^1_0(\Om)$ and $\nabla H^1_0(\Om)$, viewed as a subset of $H_0(\crl,\Om)$. The gradients of any orthonormal basis of $H^1_0(\Om)$ thus define an orthonormal basis of $\nabla H^1_0(\Om) \subset \Hcz$. 

From classical spectral theory, the following proposition defines an orthonormal basis for $\nabla H^1_0(\Om)$ in terms of the homogeneous-Dirichlet eigenvectors of $-\Delta$ in $\Om$. 
\begin{proposition}\label{propNabla}
  Let $\Om$ be a bounded, simply connected, and Lipschitz domain in $\mathbb{R}^n$, where $n = 2$ or $3$. There exists an orthonormal basis of $H^1(\Om)$, consisting of non-zero homogeneous-Dirichlet eigenvectors of $-\Delta$ in $\Om$, $(\phi_k)_{k\in\I}$, for a countable index set $\I$.  Then, the sequence $(\bphi_k)_{k\in\I}$, defined as \[\bphi_k := \frac{\nabla \phi_k}{\| \nabla \phi_k \|_{[L^2(\Om)]^n}},\] 
  forms an orthonormal basis for $\nabla H^1_0(\Om) \subset \Hcz$. 
\end{proposition}
Appendix~\ref{sec:appA} contains the proof of Proposition~\ref{propNabla}.

\subsection{A set of basis functions for $X_0(\Om)$}

For constructing the eigenbasis for the space $X_0(\Om)$, we resume the prior discussion on the eigenvectors of the operator $(1+\crl \text{-} \crl)$. We recall that  $X_0(\Om)$ is compactly embedded into $L^2(\Om)$ (see \cite[Theorem 4.23]{kirsch2015mathematical}), which ensures that the $(1+\crl \text{-} \crl)$ operator admits a compact and self-adjoint inverse when restricted to $X_0(\Om)$. Then, the eigenbasis of the operator forms a complete and orthonormal basis of the space $X_0(\Om)$. Moreover, when restricted to divergence-free vector fields, the $(1+\crl \text{-} \crl)$ operator reduces to $(1-\Delta)$, which is suggestive of the fact that one may construct eigenvectors of $(1+\crl \text{-} \crl)$ via eigenvectors of the negative Laplacian, which we perform in the following propositions. Here, we divide the construction of this eigenbasis into two scenarios based on the dimensionality of $\Om$.
\begin{proposition}\label{propX2}
  Let $\Om$ be a bounded, simply connected, and Lipschitz domain in $\mathbb{R}^2$. There exists an orthonormal basis of $H^1(\Om)$, consisting of non-constant homogeneous-Neumann eigenvectors of $-\Delta$ in $\Om$, $(\bpsi_k)_{k\in\I}$, for a countable index set $\I$.  Then, the sequence $(\bpsi_k)_{k\in\I}$, defined as \[\bpsi_k := \frac{\crl^*(\phi_k)}{\| \crl^*(\phi_k) \|_{\Hc}},\] where $\crl^*$ is as defined in \eqref{curl2}, forms an orthonormal basis for $X_0(\Om) \subset \Hcz$. 
\end{proposition}
The proof of Proposition~\ref{propX2} is detailed in Appendix~\ref{sec:appA}.

Finding an orthonormal basis for the space $X_0(\Om)$ in 3D is more complex than in the 2D case. Here, we restrict the construction of basis functions for the space $X_0(\Om)$ to Cartesian product domains, i.e., domains defined as the Cartesian product of a simply connected domain in $\RR^2$ and a closed interval in $\RR$. In this particular case, we can construct eigenvectors of the $(1+\crl \text{-} \crl)$ operator in a similar fashion, according to differential operators acting upon scalar-valued eigenfunctions of the Laplacian with appropriate boundary conditions. 

In such geometries, the basis functions of $X_0(\Om)$ in three dimensions come in two distinct modes, called TM (Transverse Magnetic) and TE (Transverse Electric) modes \cite{kirsch2015mathematical}. The construction of such an eigenbasis was developed in \cite{costabel2019maxwell}. There, the authors demonstrate that the TE and TM modes defined below constitute a complete and orthonormal basis of $X_0(\Om)$. 

Without loss of generality, we state that a Cartesian product domain $\Om\in \RR^3$ is of the form $\Om^* \times I$, with $\Om^* \subset \RR^2$ simply connected and $I \subset \RR$ an interval. We refer to the coordinate direction corresponding to the interval $I$ as the distinguished direction. For this specific definition of $\Om$, the $z$-coordinate is the distinguished axis of the Cartesian product domain, but it is worth noting that the following derivation is direction-independent.

Considering the scope of this work, we refer to \cite{costabel2019maxwell} for a more detailed explanation and in Appendix~\ref{sec:appA} we give further details of the properties of the basis functions for $X_0(\Omega)$ in 3D. 

\begin{proposition}\label{propX3}
  Let $\Om \subset \RR^2$ be a Cartesian product domain such that $\Om = \Om^* \times I$, where $\Om^* \subset \mathbb{R}^2$ is a simply connected domain, and $I\subset \RR$ is an interval. Given two sets of indices $\I_1$ and $\I_2$, consider the following sets of functions:
  \begin{itemize}
  	\item The non-zero functions $(p_{k})_{k\in\I_1}$ forming a complete set of eigenvectors of $-\Delta$ in $\Om$  with homogeneous-Dirichlet boundary conditions on $\overline{\Om^*} \times \partial I$ and homogeneous-Neumann boundary conditions on $\overline{I} \times \partial \Om^*$.
  	\item  The non-zero functions $(q_{h})_{h\in\I_2}$ forming a complete set of eigenvectors of $-\Delta$ in $\Om$ with homogeneous-Neumann boundary conditions on $\overline{\Om^*} \times \partial I$ and homogeneous-Dirichlet boundary conditions on $\overline{I} \times \partial \Om^*$.
  \end{itemize}
  Then,  we define the sequences of vector fields $(\bpsi_{k}^{\mathrm{TM}})_{{k}\in\I_1}$ and $(\bpsi_{h}^{\mathrm{TE}})_{{h}\in\I_2} \in X_0(\Om)$, via
  \[\bpsi_{k}^{\mathrm{TM}}  := \frac{ \crl(p_{k} \te)}{\| \crl(p_{k} \te) \|_{\Hc}},  \text{ and } \quad \bpsi_{h}^{\mathrm{TE}}  := \frac{ \crl(\crl(q_{h} \te)) }{\| \crl(\crl(q_{h} \te)) \|_{\Hc}},\] 
  for each ${k}\in\I_1$ and ${h}\in\I_2$, and where $\te$ is the unit vector in the distinguished direction. Then, the union of the two sequences $(\bpsi_{k}^{\mathrm{TM}})_{{k}\in\I_1}$ and $(\bpsi_{h}^{\mathrm{TE}})_{{h}\in\I_2}$ forms an orthonormal basis for $X_0(\Om) \subset \Hcz$. 
\end{proposition}

 Note that the building of basis functions for the space $X_0(\Om)$ is valid, for instance, in three-dimensional rectangular, cubic, or cylindrical domains, provided the eigenbasis of the Laplacian in $\Om^*$. The calculations are similar in all of these cases. In Tables~\ref{tab2} and \ref{tab3} we show the resulting basis functions for the spaces $\nabla H_0^1(\Om)$ and $X_0(\Om)$ on $n$-dimensional cubes $\Om = [0,\pi]^n$ for $n=2,3$. For clarity, we detail the construction of the eigenbasis on $n$-dimensional cubes in Appendix~\ref{sec:appB}. An extension to more general rectangular domains with distinct side lengths is trivial, requiring only straightforward yet tedious calculations and it is therefore omitted. 

\begin{table}[htpb]
	\begin{center}
		\renewcommand{\arraystretch}{1.5}
		{\begin{tabular}{|c|l|c|}
				\hline
				- & Basis functions & Index $k = (k_1,k_2)$ \\
				\hline
				$\nabla H_0^1(\Om)$ &  $\bphi_k = \frac{2}{\pi|k|}\begin{pmatrix}k_1\cos(k_1x)\sin(k_2y)\\ k_2\sin(k_1x)\cos(k_2y) \end{pmatrix} $ & $k_\ii\in \ZZ_{\geq 0}$   \\
				\hline
				$X_0(\Om)$ & $\bpsi_k = \frac{c_k}{\sqrt{|k|^4+|k|^2}}\begin{pmatrix} k_2 \cos(k_1x)\sin(k_2y)\\ -k_1 \sin(k_1x)\cos(k_2y) \end{pmatrix}$ & $k_\ii\in \ZZ_{\geq 0}$ {\footnotesize with $k_{1} = 0$ xor $k_{2} = 0$} \\
				\hline
		\end{tabular}}
	\end{center}
	\vspace{-12pt}
	\caption{The basis functions for the spaces $\nabla H_0^1(\Om)$ and $X_0(\Om)$ in 2D with $\Om = [0,\pi]^2$. Here, $c_k = \frac{2}{\pi} $ if $k_\ii>0$ and $c_k = \frac{\sqrt{2}}{\pi} $ if $k_1 = 0$ or $k_2 = 0$.}
	\label{tab2}
\end{table}

\begin{table}[htpb]
  \begin{center}
    \renewcommand{\arraystretch}{1.5}
    {\begin{tabular}{|c|l|M{0.2\textwidth}|}
        \hline
        - & Basis functions & Index $k = (k_1,k_2)$ \\
        \hline
        $\nabla H_0^1(\Om)$ &  $\bphi_k = \frac{2\sqrt{2}}{ \pi^{3/2}|k|}\begin{pmatrix}
          k_1\cos(k_1x)\sin(k_2y)\sin(k_3z)\\ k_2\sin(k_1x)\cos(k_2y)\sin(k_3z)\\ k_3\sin(k_1x)\sin(k_2y)\cos(k_3z) \end{pmatrix} $ & $k_\ii\in \ZZ_{\geq 0}$   \\
        \hline
        $X_0(\Om)$ & {$\bpsi_k^{\mathrm{TM}} = \frac{2\sqrt{2}}{\pi^{3/2}\sqrt{(1+|k|^2)(k_2^2+k_3^2)}}\begin{pmatrix}
            0 \\ -k_3\sin(k_1x)\cos(k_2y)\sin(k_3z)\\ k_2\sin(k_1x)\sin(k_2y)\cos(k_3z)
          \end{pmatrix}$ } & $k_\ii\in \ZZ_{\geq 0}$ {\footnotesize with $k_1>0$ and $k_{2} = 0$ xor $k_{3} = 0$} \\
         & {$\bpsi_k^{\mathrm{TE}}  = \frac{1}{\sqrt{c'_k(1+|k|^2)|k|^2\left(|k|^2-k_1^2\right)}}
           \begin{pmatrix}
             (k_2^2+k_3^2)\cos(k_1x)\sin(k_2y)\sin(k_3z)\\
             -k_1k_2\sin(k_1x)\cos(k_2y)\sin(k_3z)\\
             -k_1k_3\sin(k_1x)\sin(k_2y)\cos(k_3z)
           \end{pmatrix}$} & $k_\ii\in \ZZ_{\geq 0}$ {\footnotesize with $k_{2}>0$ and $k_{3}>0$}\\
        \hline
    \end{tabular}}
  \end{center}
  \vspace{-12pt}
  \caption{The basis functions for the spaces $\nabla H_0^1(\Om)$ and $X_0(\Om)$ in 3D with $\Om = [0,\pi]^3$. Here, $c'_k = \frac{\pi^3}{4}$ if $k_1 = 0$ and $c'_k = \frac{\pi^3}{8}$ if $k_1\neq0$.}
  \label{tab3}
\end{table}

\section{The DFR method}\label{sec:dfrNN}

In this section we outline the structure of the Neural Networks that we will employ and the fundamental principle underlying the DFR method, which is the construction of a discretized and computable loss function.

\subsection{Neural Networks}

A neural network is a mathematical model comprising multiple compositions of simple functions, called layers. Specifically, a NN is a non-linear function $\N(\tx; \tW, \tb)$ parametrized by a set of weights $\tW$ and biases $\tb$, with $\tx \in \RR^d$ being the input vector. We restrict our attention to fully-connected feed-forward NNs. In this architecture, the weights $\tW$ are represented as a collection of dense matrices ${ \tW_{1}, \tW_{2},..., \tW_{L} }$, where $\tW_\jj \in \RR^{d_\jj \times d_{\jj-1}}$ is the weight matrix for the layer $\jj$ and $d_\jj$ denotes the number of nodes on each layer. And the set of biases $\tb$ is a collection of vectors ${ \tb_{1}, \tb_{2},..., \tb_{L} }$, where $\tb_\jj \in \RR^{d_\jj}$ is the bias vector for the layer $\jj$.

The output of each layer of the network is computed as
\[\tA_\jj = \sigma_\jj(\tW_\jj\tA_{\jj-1} + \tb_\jj ) \]
where $\sigma_\jj$ is a nonlinear activation function, $\tA_{0} = \tx$ and $\sigma_L$ is the identity function. 
The weights and biases are obtained via a gradient-based optimization algorithm applied to a loss function, whose gradients are efficiently calculated via backpropagation. This process adjusts the parameters to minimize a specified loss function, which in our case is a discretized and computable approximation to the dual norm of the PDE residual, as defined in the following section.

We impose homogeneous Dirichlet-type boundary conditions on our candidate solutions by using a cutoff function $\xi: \overline{\Omega} \to \RR^{d\times d}$. Specifically, we define the approximation of the solution of \eqref{eq:maxwell} as $\tE = \xi\tilde{\tE}$, where $\tilde{\tE}$ is the output of the fully-connected feed-forward neural network $\N(\tx; \tW, \tb)$ and $\xi$ is smooth and non-trainable. The function $\xi$ gives matrices that are positive definite in $\Om$, and enforces the tangential operator to be zero on the boundary, i.e., the constraint $\gamma_t(\tE) = 0$ on $\Gamma$. 

\subsection{The discretized loss}

In this section, we construct a discretized and computable loss function $\LL$ as an approximation of the dual norm of the residual \eqref{eq:parseMax}. For the sake of simplicity and exposition, we restrict the details in this section to squares, i.e., $\Om = [0,\pi]^2$ using the basis functions in Table~\ref{tab2}. The upcoming calculations are similar for other Cartesian product domains in 2D and 3D. 

Notice that calculating the integrals in $\langle  \R(\tE), \Phi_k \rangle_{H'\times H}$ may be costly since it requires the evaluation of many functions and their derivatives. However, we underline that the basis functions for the space $\Hcz$ are specified in terms of Laplacian's eigenbasis and thus both the basis functions and their derivatives are described in terms of sines and cosines. This feature motivates the use of efficient methods for approximating such integrals, which are naturally based on the Fast Fourier Transform (FFT). In this regard, applying Discrete Sine/Cosine transformations (DST/DCT) is a well-known strategy that significantly reduces computational complexity when applying a $N$ point midpoint rule for $N$ basis functions. Specifically, using DST/DCT as a quadrature rule reduces the amount of required calculations from $O(N^2)$ to $O(N\log(N))$. In our case, the Discrete Sine/Cosine transforms appear naturally when one applies the mid-point integration rule to the integrals appearing in the residual. Here, we use the type II Sine/Cosine transforms defined in \cite[Section 4.2]{britanak2010discrete}. Each transform is represented by an $N \times N$ matrix as
\begin{equation*}\label{eq:mid}
	\begin{aligned}
		(S^{II}_{N})_{\ii\jj} & := \sqrt{\frac{2}{N}}\sigma_\ii \sin\left(\frac{\pi}{N}\left(\jj+\frac{1}{2} \right)(\ii+1) \right) \text{ with }
		\sigma_\ii = \left\{\begin{array}{c c}
			\frac{1}{\sqrt{2}} & \ii = N-1,\\
			1 & \ii\neq N-1,
		\end{array}\right. \\
		(C^{II}_{N})_{\ii\jj} & := \sqrt{\frac{2}{N}}\sigma'_\ii \cos\left(\frac{\pi}{N}\left( \jj+\frac{1}{2} \right)\ii \right) \text{ with }
		\sigma'_\ii = \left\{\begin{array}{c c}
			\frac{1}{\sqrt{2}} & \ii = 0,\\
			1 & \ii \neq 0,
		\end{array}\right.
	\end{aligned}
\end{equation*}
where $\ii,\jj = 0,\dots,N-1$. 

Given cutoff frequencies (Fourier modes) $N>0$ and $M>0$, we aim to approximate $\langle \R(\tE), \Phi_k \rangle_{H'\times H}$, where $k = (k_1,k_2)$ is such that $0\leq k_1\leq N$ and $0\leq k_2\leq M$. Moreover, the indices $k$ satisfy the conditions in Table~\ref{tab2}, and we let $\I$ denote the set of all appropriate indices. For each $k \in \I$, the basis functions $\Phi_k$ are $\Phi_k = \bphi_k+ \bpsi_k$, with $\bphi_k$ and $\bpsi_k$ as defined in Table~\ref{tab2}. Next, we consider the mid-point integration rule, that is, 
\begin{equation}\label{eq:mid}
	\int_{\Om} f \, d\tx \approx \frac{\pi^2}{NM} \sum_{\ii = 0}^{N-1}\sum_{\jj = 0}^{M-1} f(x_\ii,y_\jj), 
\end{equation}
where the integration points are $x_\ii = \frac{2\ii+1}{2N}\pi$ and $y_\jj = \frac{2\jj+1}{2M}\pi$. Applying \eqref{eq:mid} in \eqref{eq:maxshort}, we approximate the integrals appearing in the residual as 
\begin{equation}
	\begin{aligned}
		\int_\Om \mu^{-1} \crl(\tE) \cdot \crl(\Phi_k) \, d\tx &\approx \sum_{\ii = 0}^{N-1}\sum_{\jj = 0}^{M-1} \frac{\pi^2}{\sqrt{NM}} \mu^{-1} \crl( \tE)\!\!\left(x_\ii, y_\jj \right) \alpha_k \left( C^{II}_{N}\right)_{k_1\!-\!1 \ii} \left( C^{II}_{M}\right)_{k_2\!-\!1 \jj} \\
		& = : \mathcal{R}_{1k}(\mu^{-1} \crl(\tE))\\
		\int_\Om \omega^2 \epsilon \tE \cdot \Phi_k \, d\tx &\approx \sum_{\ii = 0}^{N-1}\sum_{\jj = 0}^{M-1} \frac{\pi^2}{\sqrt{NM}} \omega^2\epsilon \, \tE\!\left( x_\ii, y_\jj \right)\cdot\bal'_k \mathbf{C}^{k}_{\ii\jj}\mathbf{S}^{k}_{\ii\jj}\\
		& = : \mathcal{R}_{2k}(\omega^2\epsilon \, \tE)\\
		\int_\Om  i \omega \tJ \cdot \Phi_k \, d\tx &\approx \sum_{\ii = 0}^{N-1}\sum_{\jj = 0}^{M-1} \frac{\pi^2}{\sqrt{NM}} i \omega \tJ \!\left( x_\ii, y_\jj \right)\cdot\bal'_k \mathbf{C}^{k}_{\ii\jj}\mathbf{S}^{k}_{\ii\jj}\\
		& = \mathcal{R}_{2k}(i \omega \tJ)
	\end{aligned}
\end{equation}
where $\alpha_k = \frac{c_k (k_1k_2-k_1^2)}{\sqrt{|k|^4+|k|^2}}$ and $ \bal'_k = \left( \frac{2k_1}{\pi|k|} + \frac{c_k k_2}{\sqrt{|k|^4+|k|^2}}, \frac{2k_2}{\pi|k|} - \frac{c_k k_1}{\sqrt{|k|^4+|k|^2}} \right)^t\!\! \mathbb{I}$, with $\mathbb{I}$ being the identity matrix. The matrices $\mathbf{C}^{k}_{\ii\jj} $ and $\mathbf{S}^{k}_{\ii\jj} $ contain the cosine and sine transformations as follows 
\begin{equation*}
	\mathbf{C}^{k}_{\ii\jj} = \begin{pmatrix}
		\left( C^{II}_{N}\right)_{k_1\!-\!1 \ii} & 0\\
		0 & ( C^{II}_{M})_{k_2\!-\!1 \jj}
	\end{pmatrix} \, \text{ and } \, \mathbf{S}^{k}_{\ii\jj} = \begin{pmatrix}
		( S^{II}_{M})_{k_2\!-\!1 \jj} & 0 \\
		0 & \left( S^{II}_{N}\right)_{k_1\!-\!1 \ii} 
	\end{pmatrix}.
\end{equation*}
Notice that, to evaluate the $\crl$ of the candidate solution $\tE$ one needs to use automatic differentiation as described in \cite{baydin2018automatic}. Finally, our discretized loss is 
\begin{equation}\label{eq:discloss}
	 \LL(\tE)  := \sqrt{\sum_{k\in \I} |\mathcal{R}_{1k}(\mu^{-1} \crl(\tE)) - \mathcal{R}_{2k}(\omega^2\epsilon \, \tE) - \mathcal{R}_{2k}(i \omega \tJ)|}.
\end{equation}

We point out that it is also possible to employ a different number of integration points and basis functions (modes), but we omit the details here. Notice that the use of less integration points that Fourier modes is expected to reduce integration errors.

\section{Numerical experiments}\label{sec:numerical}
In this section, we present four numerical experiments that illustrate the main capabilities and some of the limitations of the DFR method. We use \textit{Tensorflow 2.8} and implement a feed-forward fully-connected NN that consists of five hidden layers, each with $20$ neurons with a $\tanh$ activation function. In our implementation, Adam serves as the optimizer.  Through the use of a \textit{Callback}, we allow the optimiser to dynamically modify the learning rate based on the decay of the loss and reject iteration steps that result in an increase in loss (see \cite{uriarte2022finite}). We choose an starting learning rate of $10^{-4}$. In addition, we use a validation set and compare the training and validation losses every iteration. This method aids in detecting overfitting, as mentioned in \cite{rivera2022quadrature}.

\subsection{Case 1. Smooth solution in 2D}\label{sec:case1}
Let $\Om =[0,\pi]^2$ and consider the variational form: find $\tE \in \Hcz$ satisfying
\begin{equation}\label{eq:ex1}
	\int_\Om  \crl(\tE) \cdot \crl(\bphi)+ \tE \cdot \bphi \, d\tx = \int_\Om \tilde{\tJ} \cdot \bphi \, d\tx \qquad \forall \bphi \in \Hcz.
\end{equation}
Here, $\tilde{\tJ}$ is chosen such that the exact solution is $\tE^*(x,y) = (xy(y-\pi),xy(x-\pi))^t$. 
Notice that, the bilinear form in \eqref{eq:ex1} is precisely the inner product on $\Hc$ and thus the dual norm of the PDE residual and the $H$(curl)-norm of the error are equal, in the sense that the constants $M$ and $\gamma$ in \eqref{eq:bounds} are equal to $1$. Whilst this physically does not correspond to the standard Maxwell's equations, we include this as a test of the mathematical accuracy of the method.

 We take $N=M=100$ integration points for training the NN and $117$ integration points on each direction for validation. Moreover, we also select $100$ modes in both training and validation. 
 
\begin{figure}[htbp!]
 	\centering
 	\begin{subfigure}[b]{0.495\textwidth}
 		\centering
 		\includegraphics[width=1\textwidth]{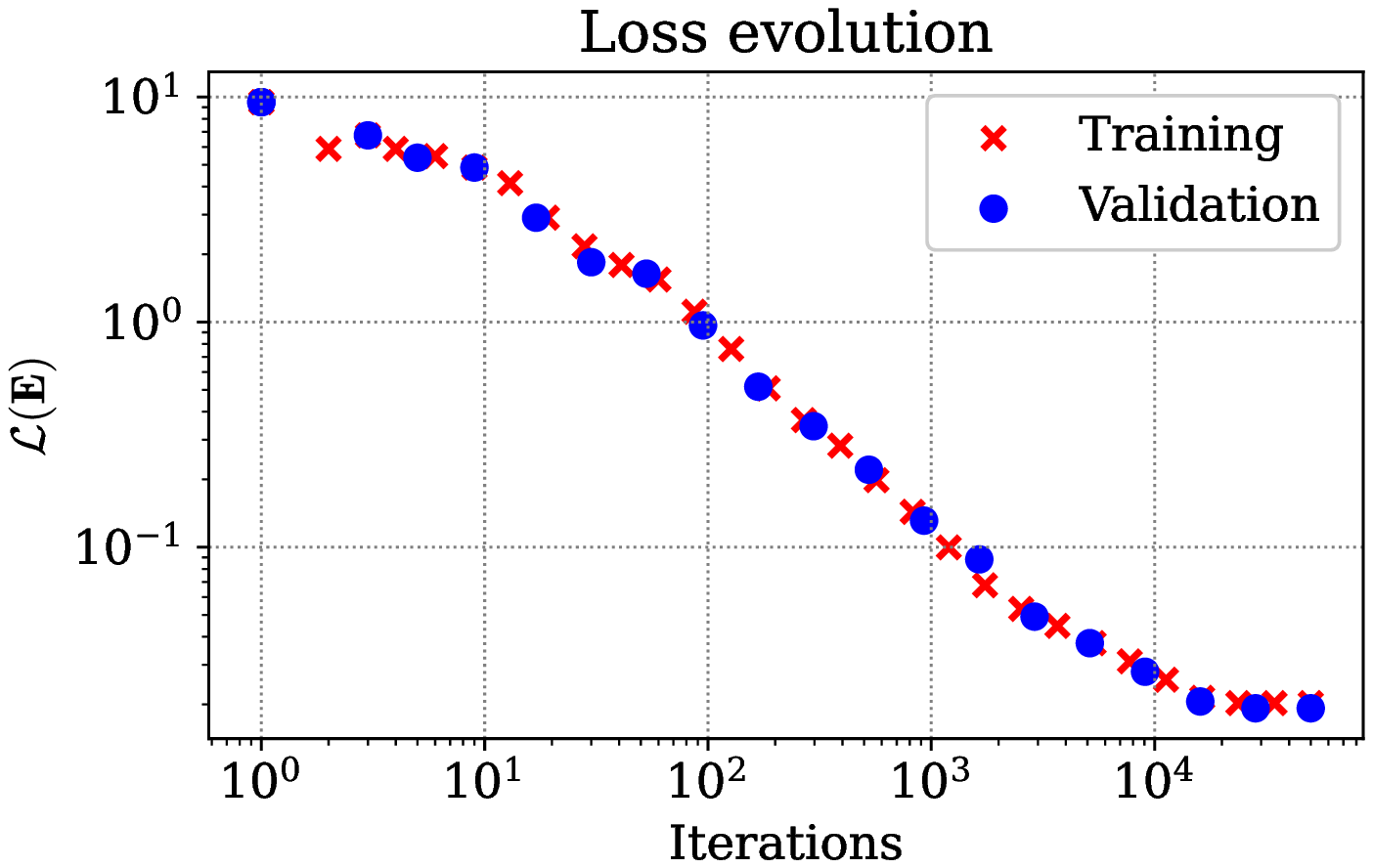}
 		\caption{\centering The evolution of the loss $\LL(\tE)$ on both the training and validation data sets.}
 		\label{fig:lossev_ex1}
 	\end{subfigure}
 	\begin{subfigure}[b]{0.495\textwidth}
 		\centering
 		\includegraphics[width=1\textwidth]{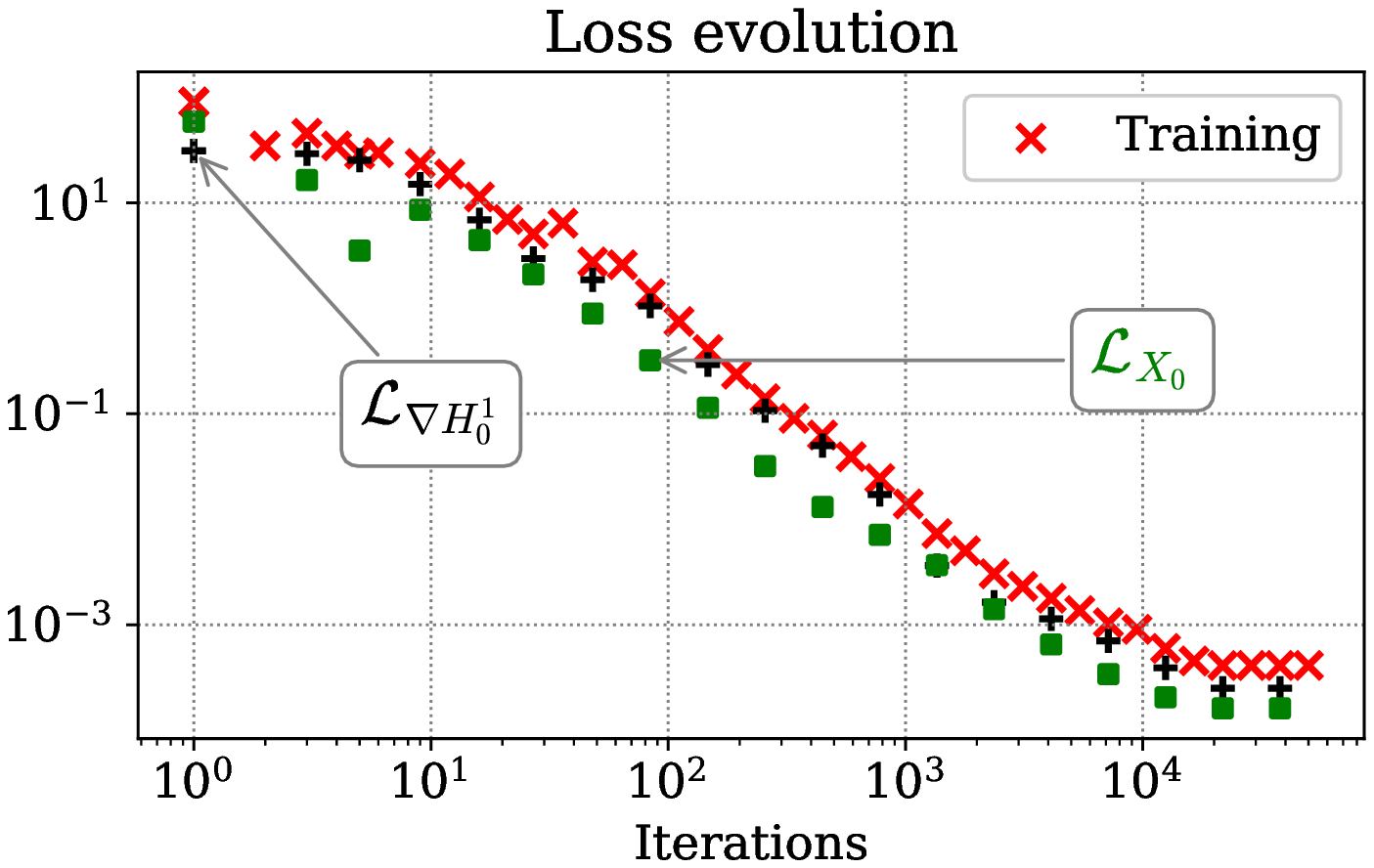}
 		\caption{ \centering The contribution of the spaces $\nabla H_0^1(\Om)$ and $X_0(\Om)$ to the loss  $\LL(\tE)^2$ on the training set.}
 		\label{fig:lossev2_ex1}
 	\end{subfigure}
 	\caption{The evolution of the loss in Case~\hyperref[sec:case1]{1}.}
 	\label{fig:P1_2}
 \end{figure}

 In Figure~\ref{fig:lossev_ex1} we show the evolution of the loss on the training and validation data sets. After $10^4$ iterations, both losses stabilize and converge to a limiting value. On the other hand, we show in Figure~\ref{fig:lossev2_ex1} the contribution of the spaces $\nabla H_0^1(\Om)$ and $X_0(\Om)$ to the squared loss on the training set. We let $\LL_{\nabla H_0^1}$ and $\LL_{X_0}$ denote the parts of $\LL(\tE)^2$ computed by using only the basis functions of the corresponding subspaces. 
 
\begin{figure}[htbp!]
 	\centering
 	\includegraphics[width=0.6\textwidth]{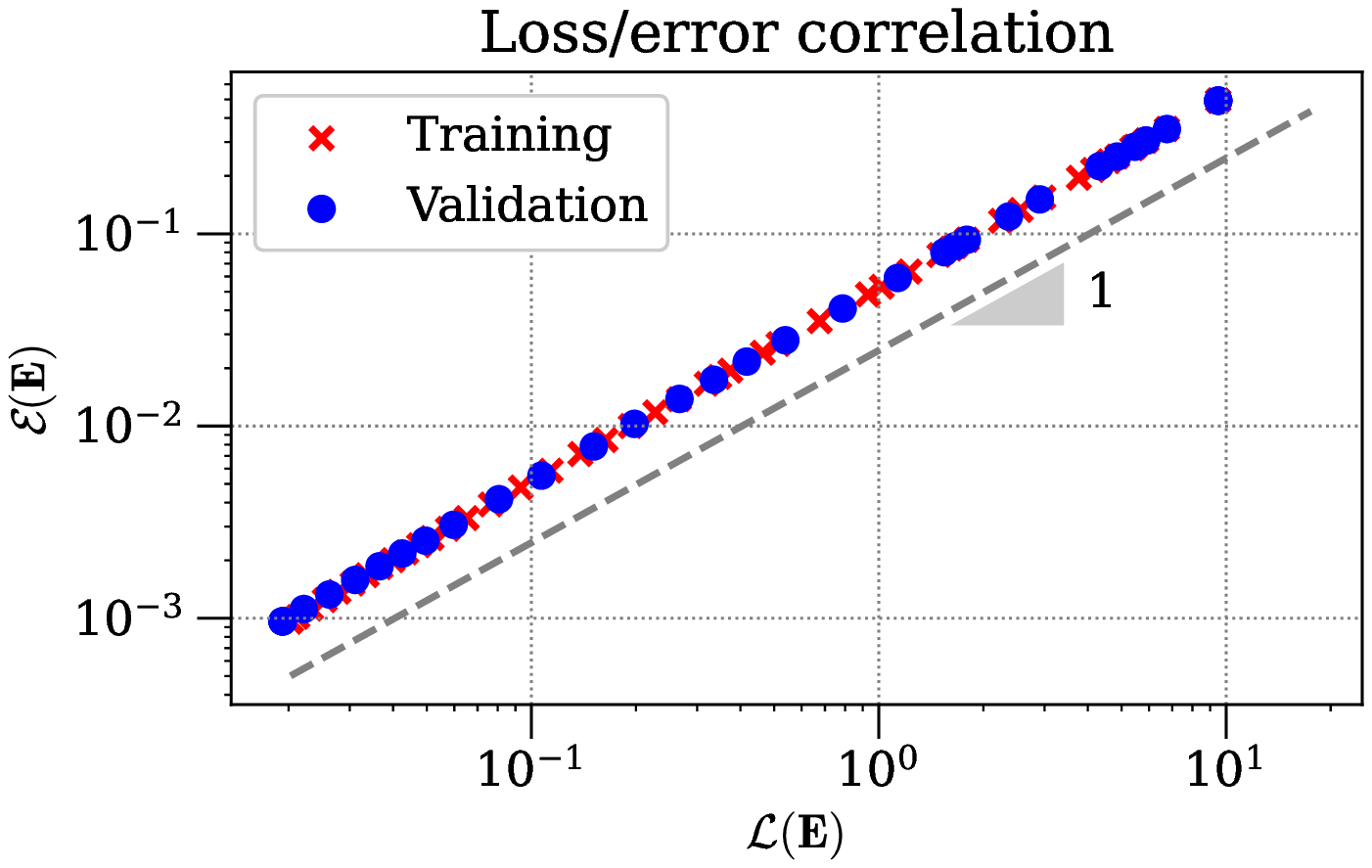}
 	\vspace{-0.6cm}
 	\caption{\centering The correlation between the loss $\LL(\tE)$ and the relative error of the solution $\mathcal{E}(\tE)$ during training and validation in Case~\hyperref[sec:case1]{1}.}
 	\label{fig:P1_3}
 \end{figure}

 In Figure~\ref{fig:P1_3} we show the relationship between the losses and the relative error of the solution during training and validation. We define the relative error in terms of the $H$(curl)-norm of the error as: \[ \mathcal{E}(\tE) := \frac{\| \tE -\tE^*\|_{\Hc} }{ \| \tE^*\|_{\Hc} } \] and we always measure this error on the validation set.  In Figure~\ref{fig:P1_3}  we include an straight line with slope one and highlight the linear relationship between the loss and the relative error, as expected in view of \eqref{eq:bounds}. We conclude that the proposed loss is an accurate approximation of the $H$(curl)-norm of the error.
 
\begin{figure}[htbp!]
	\centering
	\begin{subfigure}[b]{0.48\textwidth}
		\centering
		\includegraphics[width=1.1\textwidth]{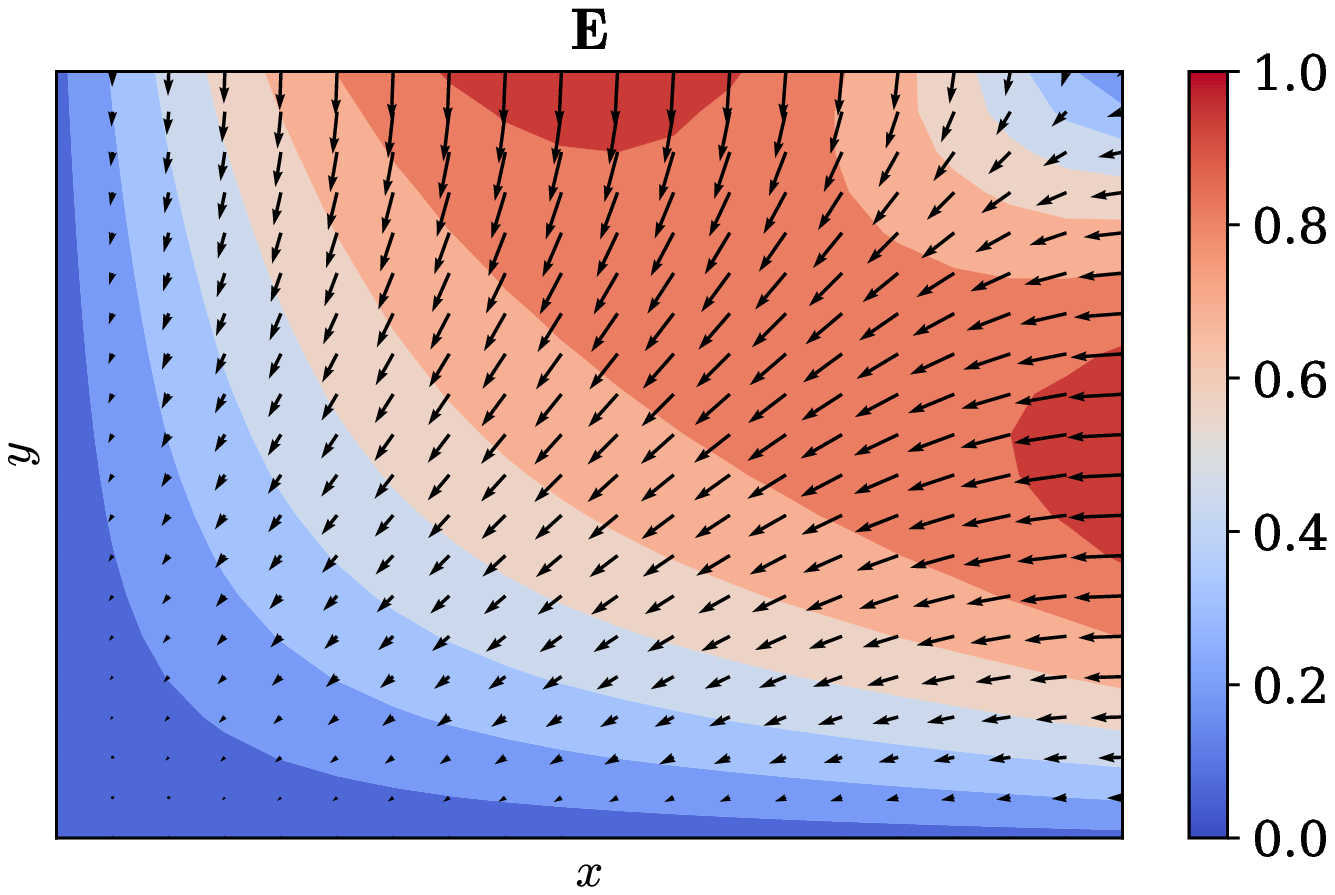}
		\caption{Approximate solution.}
		\label{fig:sol_ex1}
	\end{subfigure}
	\begin{subfigure}[b]{0.48\textwidth}
		\centering
		\includegraphics[width=1.1\textwidth]{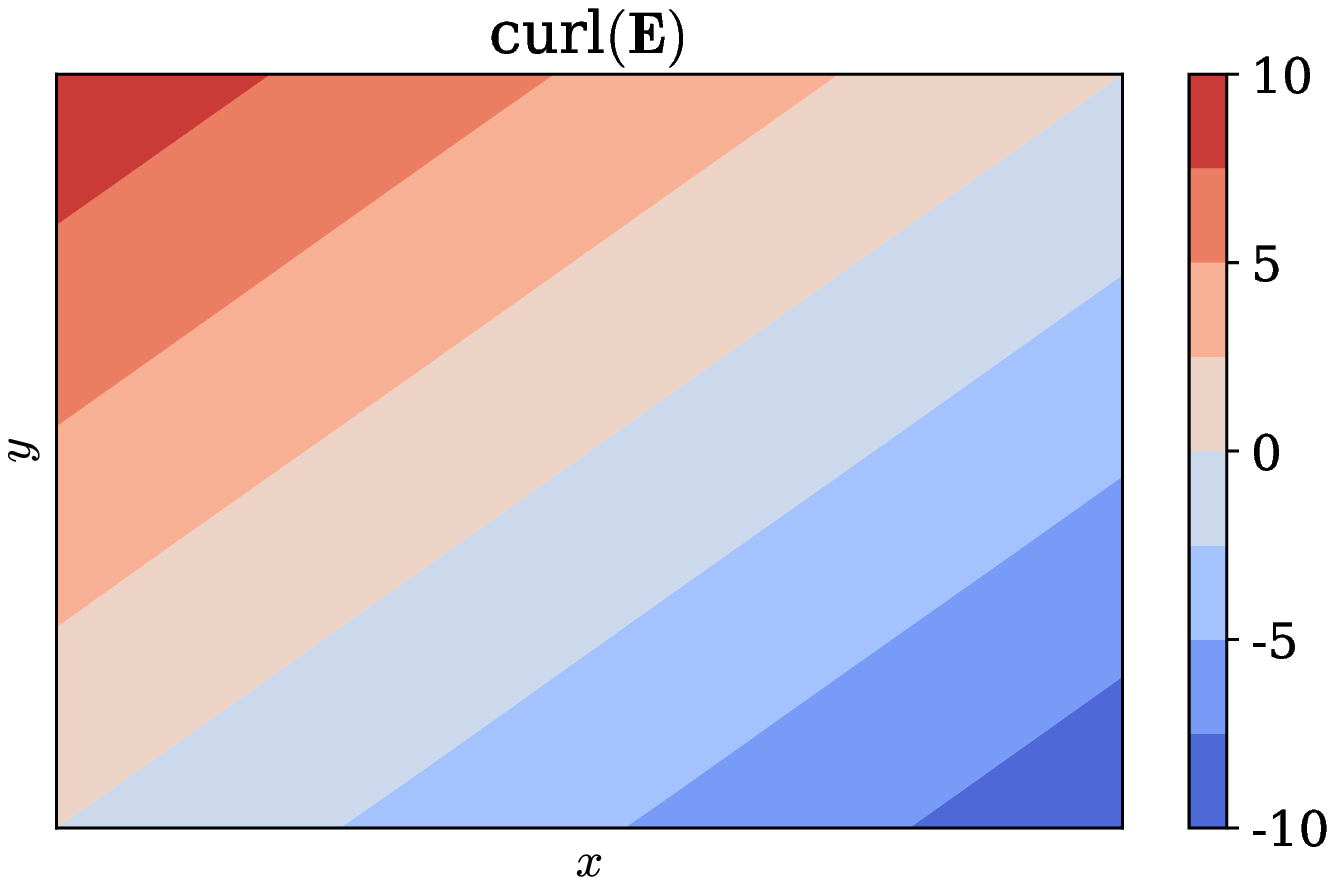}
		\caption{Curl of the approximate solution.}
		\label{fig:curl_ex1}
	\end{subfigure}
	\\
	\begin{subfigure}[b]{0.48\textwidth}
		\centering
		\includegraphics[width=1.1\textwidth]{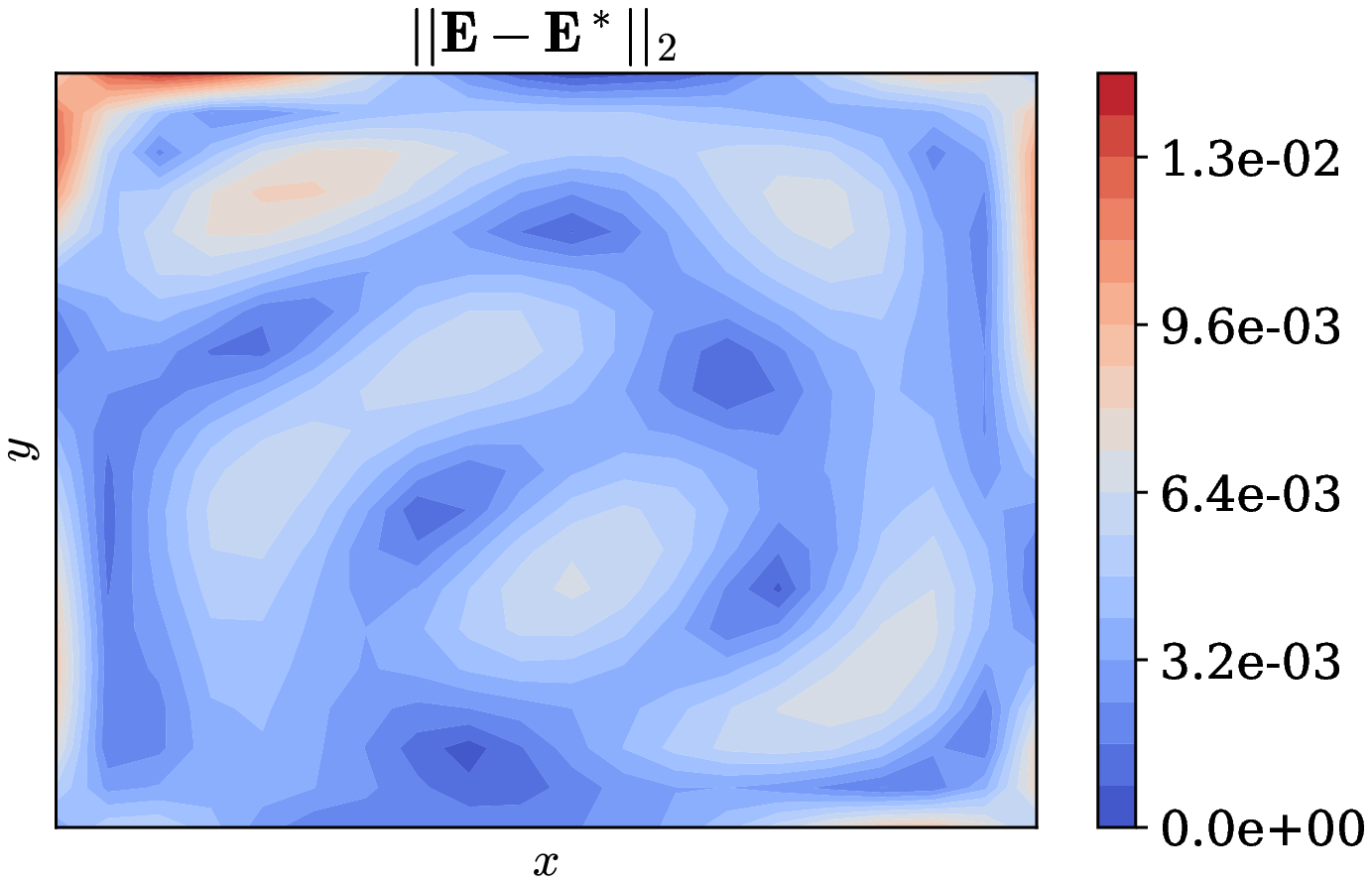}
		\caption{Error in the solution.}
		\label{fig:errorsol_ex1}
	\end{subfigure}
	\begin{subfigure}[b]{0.48\textwidth}
		\centering
		\includegraphics[width=1.1\textwidth]{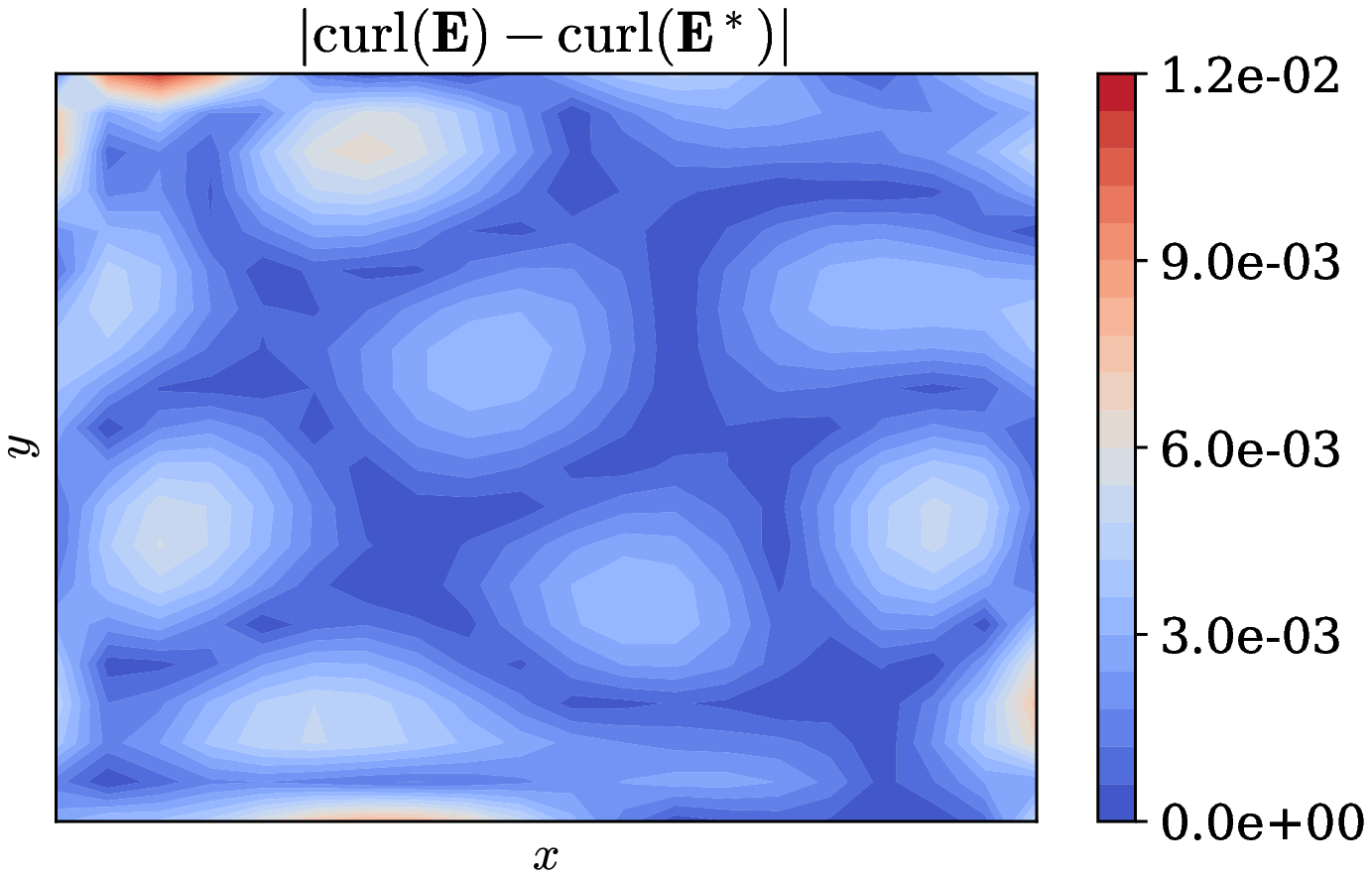}
		\caption{Error in the curl of the solution.}
		\label{fig:errorcurl_ex1}
	\end{subfigure}
	\caption{Solution and errors for the model Case~\hyperref[sec:case1]{1}.}
	\label{fig:P1_1}
\end{figure}

 Figure~\ref{fig:P1_1} shows the obtained solution, the curl of the solution, and the corresponding errors calculated pointwise.
 
\subsection{Case 2. Discontinuous parameters in 2D}\label{sec:case20}

We take $\Om =[0,\pi]^2$. We define $\Omega_0=\left\{(x,y):f_+(x,y)<1\right\}$ with $f_+(x,y)=\left(x-\frac{\pi}{2}\right)^2+\left(y-\frac{\pi}{2}\right)^2$, and select $\mu$ and $\epsilon$ to be piecewise continuous, such that 
\begin{equation*}
	\begin{split}
		\mu(x,y)=&\mu_1{\bf 1}_{\Omega_0}+\mu_2(1-{\bf 1}_{\Omega_0}),\\
		\epsilon(x,y)=&\epsilon_1{\bf 1}_{\Omega_0}+\epsilon_2(1-{\bf 1}_{\Omega_0}),
	\end{split}
\end{equation*}
with $\mu_1=3,\mu_2=1$, $\epsilon_1=1$ and $\epsilon_2=3$. With these particular parameters, we will consider two PDEs, one corresponding the coercive variational form, and a particular case of Maxwell's equations. In both cases, we seek for an exact solution $\tE^*=(E_1^*, E_2^*)^t$ defined as:
 \begin{equation}\label{eq:exE}
	\begin{split}
		E_1^*(x,y)=& \left\{\begin{array}{c c} -\mu_1 k_1\left(1-f_-(x,y)\right)\left(y-\frac{\pi}{2}\right) & (x,y)\in\Omega_0\\
			-\mu_2 k_2\left(1-f_-(x,y)\right)\left(r^2-f_-(x,y)\right)\left(y-\frac{\pi}{2}\right) & \text{else}
		\end{array} \right.\\
		E_2^*(x,y) = &\left\{\begin{array}{c c} -\mu_1 k_1\left(1-f_-(x,y)\right)\left(x-\frac{\pi}{2}\right) & (x,y)\in\Omega_0\\
			-\mu_2 k_2\left(1-f_-(x,y)\right)\left(r^2-f_-(x,y)\right)\left(x-\frac{\pi}{2}\right) & \text{else}
		\end{array} \right.\\
	\end{split}
\end{equation}
where $f_-(x,y)=\left(x-\frac{\pi}{2}\right)^2-\left(y-\frac{\pi}{2}\right)^2$, $k_1=1,k_2=35$ and $r=6$. In particular, we note that the solution admits discontinuities in both the vector field itself and the curl across the $\partial \Omega_0$, and thus the solution is in $\Hc\setminus H^1(\Om)$.

\subsection{Case 2.1. Coercive variational form}\label{sec:case21}

First, we consider the variational form: find $\tE \in \Hcz$ satisfying
\begin{equation}\label{eq:case2}
	\int_\Om \mu^{-1}  \crl(\tE) \cdot \crl(\bphi)+ \epsilon \tE \cdot \bphi \, d\tx = \int_\Om \tilde{\tJ} \cdot \bphi \, d\tx \qquad \forall \bphi \in \Hcz.
\end{equation}
Here, $\tilde{\tJ}$ is chosen such that the exact solution is \eqref{eq:exE}. Notice that the exact solution does not satisfy the strong form of the equation due to interface conditions at the discontinuities.  

In this case, as the loss requires the integration of discontinuous functions, we expect more significant integration errors during training. Thus, we employ a larger number of integration points than modes to mitigate this problem. Specifically, we take $N=M=200$ integration points for training the NN and $234$ points on each direction for validation. In both cases, we use $150$ modes.

\begin{figure}[htbp!]
	\centering
	\begin{subfigure}[b]{0.495\textwidth}
		\centering
		\includegraphics[width=1\textwidth]{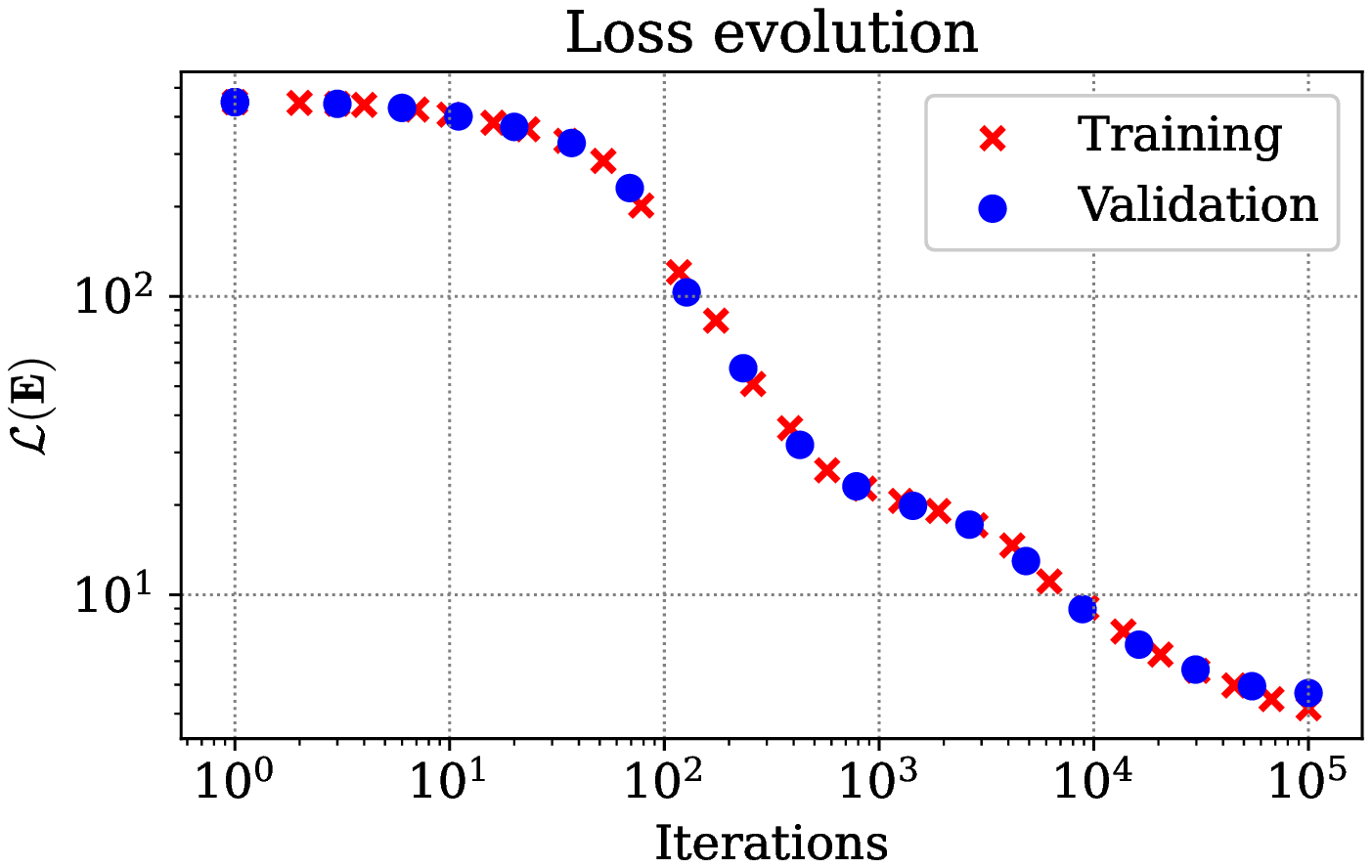}
		\caption{\centering The evolution of the loss $\LL(\tE)$ on both the training and validation data sets. }
		\label{fig:lossev_ex2}
	\end{subfigure}
	\begin{subfigure}[b]{0.495\textwidth}
		\centering
		\includegraphics[width=1\textwidth]{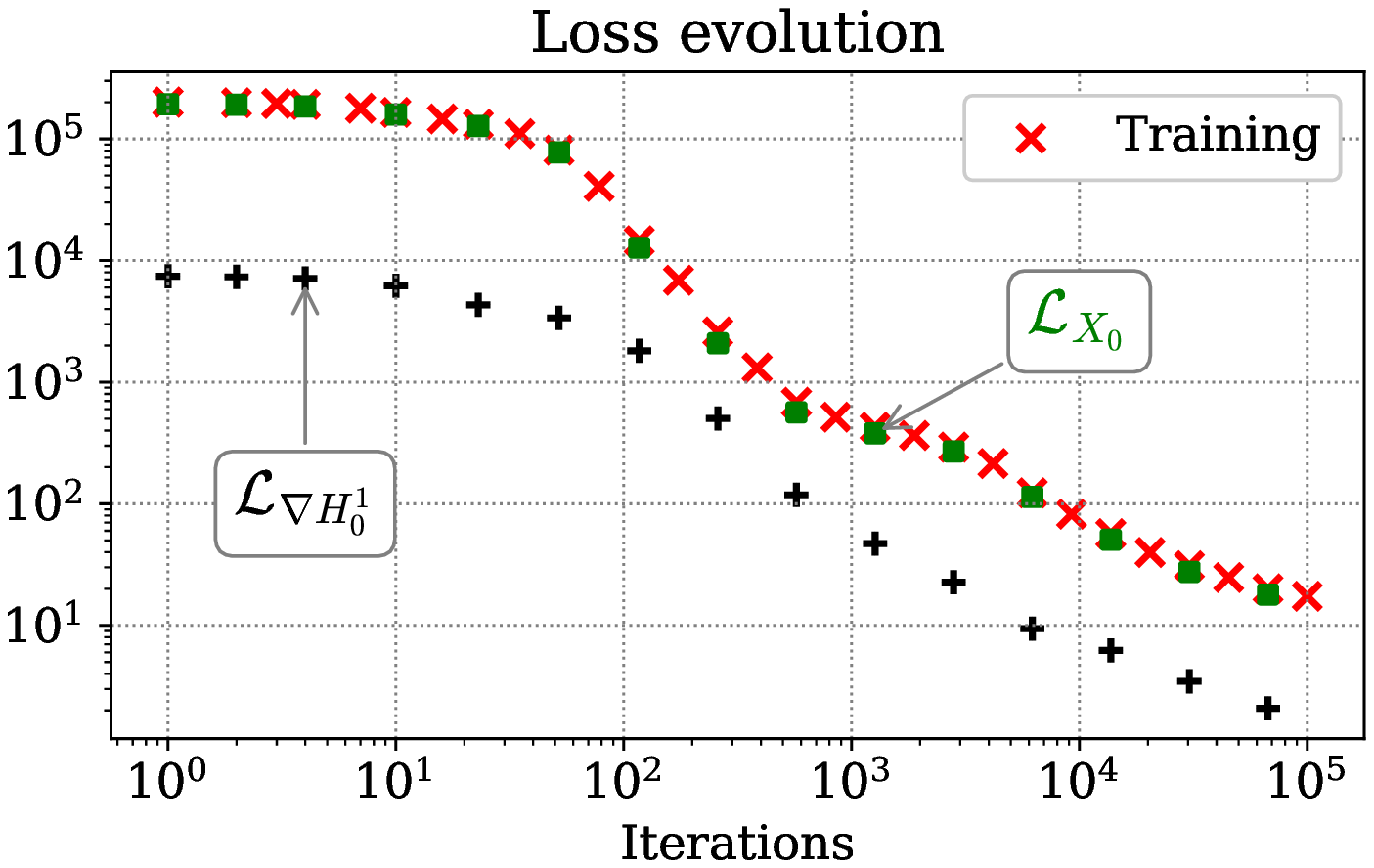}
		\caption{ \centering The contribution of the spaces $\nabla H_0^1(\Om)$ and $X_0(\Om)$ to the loss $\LL(\tE)^2$ on the training set.}
		\label{fig:lossev2_ex2}
	\end{subfigure}
	\caption{The evolution of the loss in Case~\hyperref[sec:case21]{2.1}.}
	\label{fig:P2_2}
\end{figure}

The evolution of the loss on the training and validation data sets is shown in Figure~\ref{fig:lossev_ex2}.  Figure~\ref{fig:lossev2_ex2} shows the contribution of the spaces $\nabla H_0^1(\Om)$ and $X_0(\Om)$ to the squared loss on the training set, denoted $\LL_{\nabla H_0^1}$ and $\LL_{X_0}$, respectively. The largest loss is attributed to the space ${X_0(\Omega)}$. This means that any enhancements to the basis functions in this domain will have a major impact on lowering the overall loss. Further research in this area could include parameterizing the number of basis functions related to each space separately. 

\begin{figure}[htbp!]
	\centering
	\begin{subfigure}[b]{0.48\textwidth}
		\centering
		\includegraphics[width=1.1\textwidth]{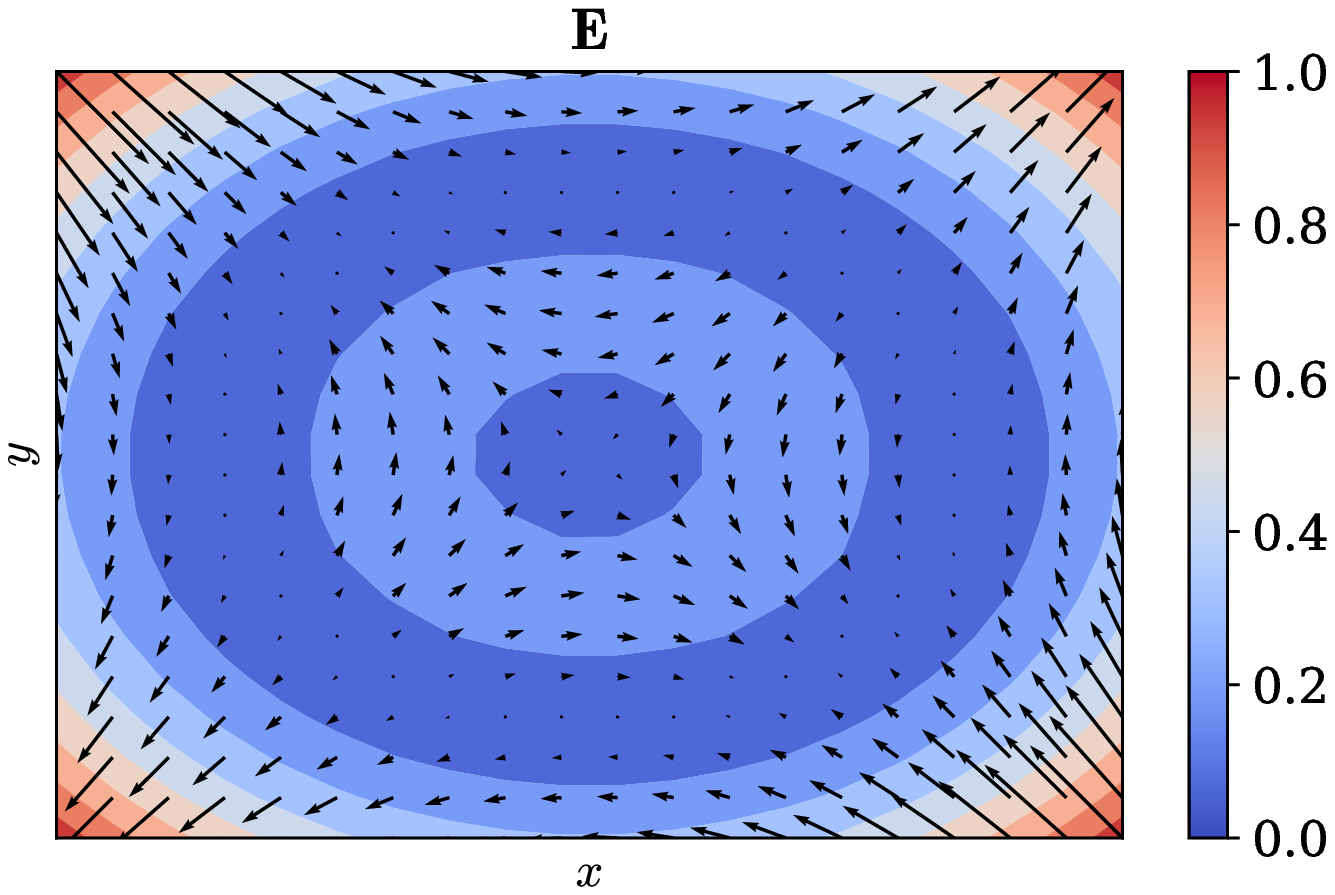}
		\caption{Approximate solution.}
		\label{fig:sol_ex2}
	\end{subfigure}
	\begin{subfigure}[b]{0.48\textwidth}
		\centering
		\includegraphics[width=1.1\textwidth]{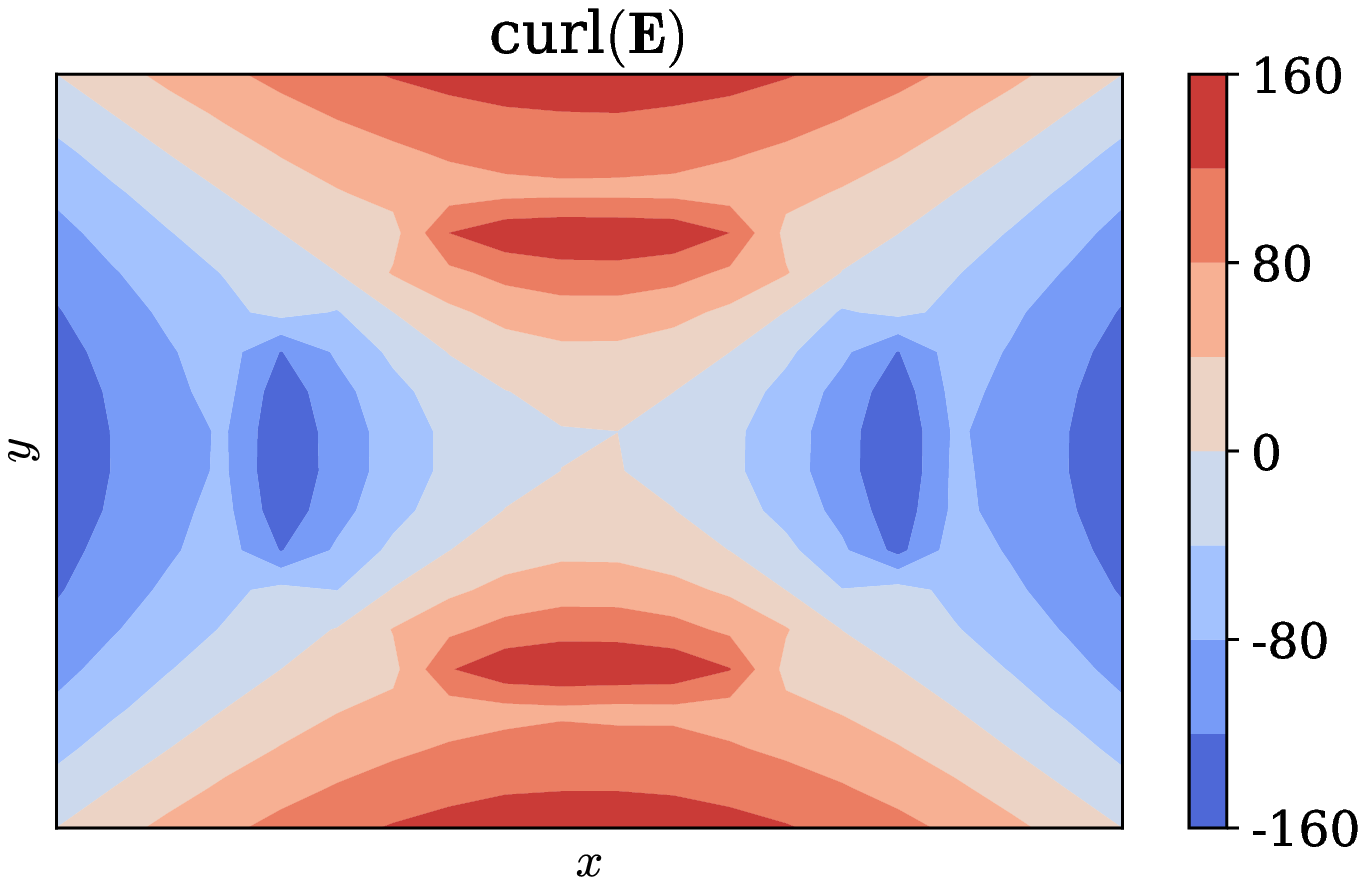}
		\caption{Curl of the approximate solution.}
		\label{fig:curl_ex2}
	\end{subfigure}
	\\
	\begin{subfigure}[b]{0.48\textwidth}
		\centering
		\includegraphics[width=1.1\textwidth]{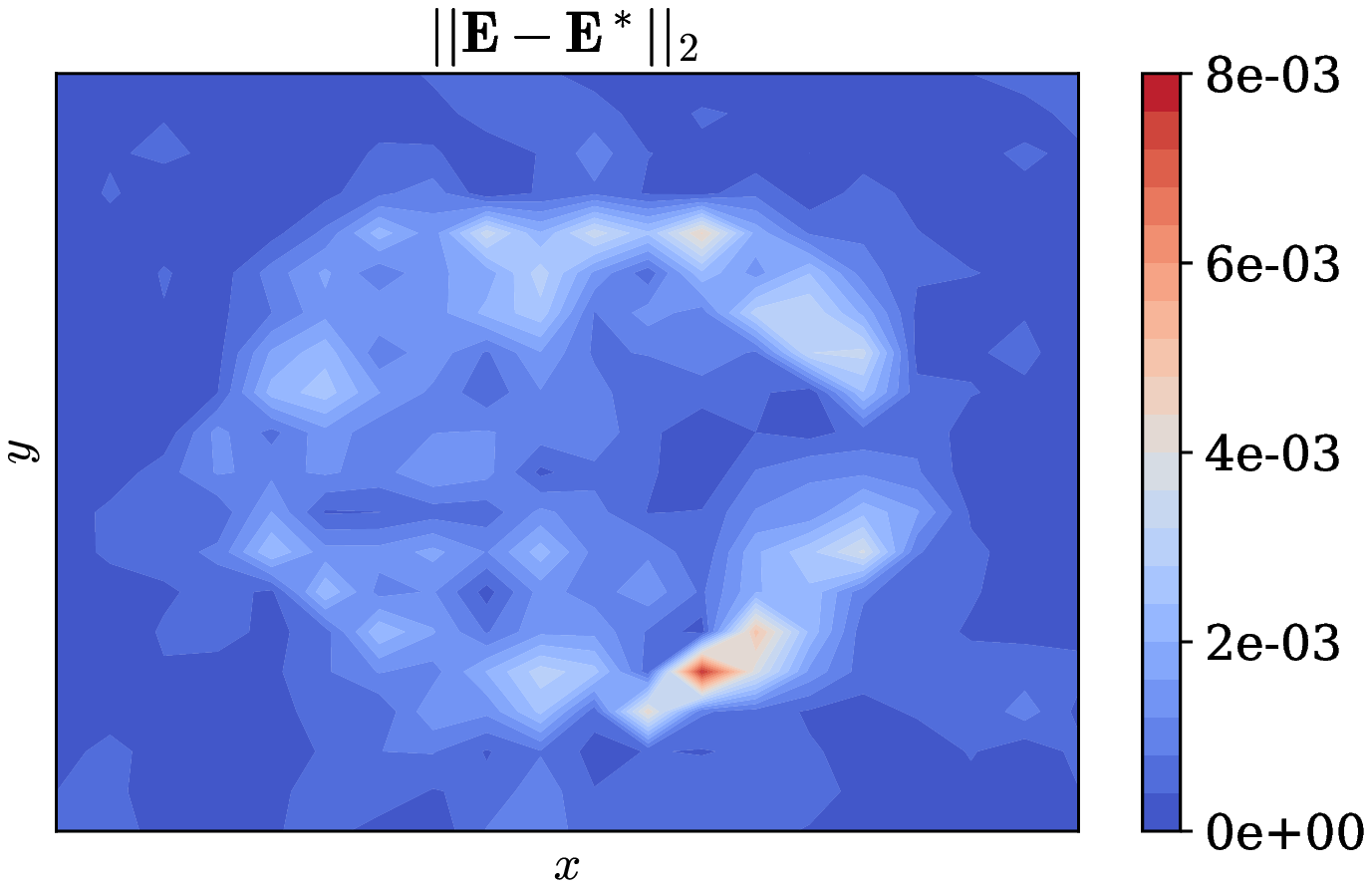}
		\caption{Error in the solution.}
		\label{fig:errorsol_ex2}
	\end{subfigure}
	\begin{subfigure}[b]{0.48\textwidth}
		\centering
		\includegraphics[width=1.1\textwidth]{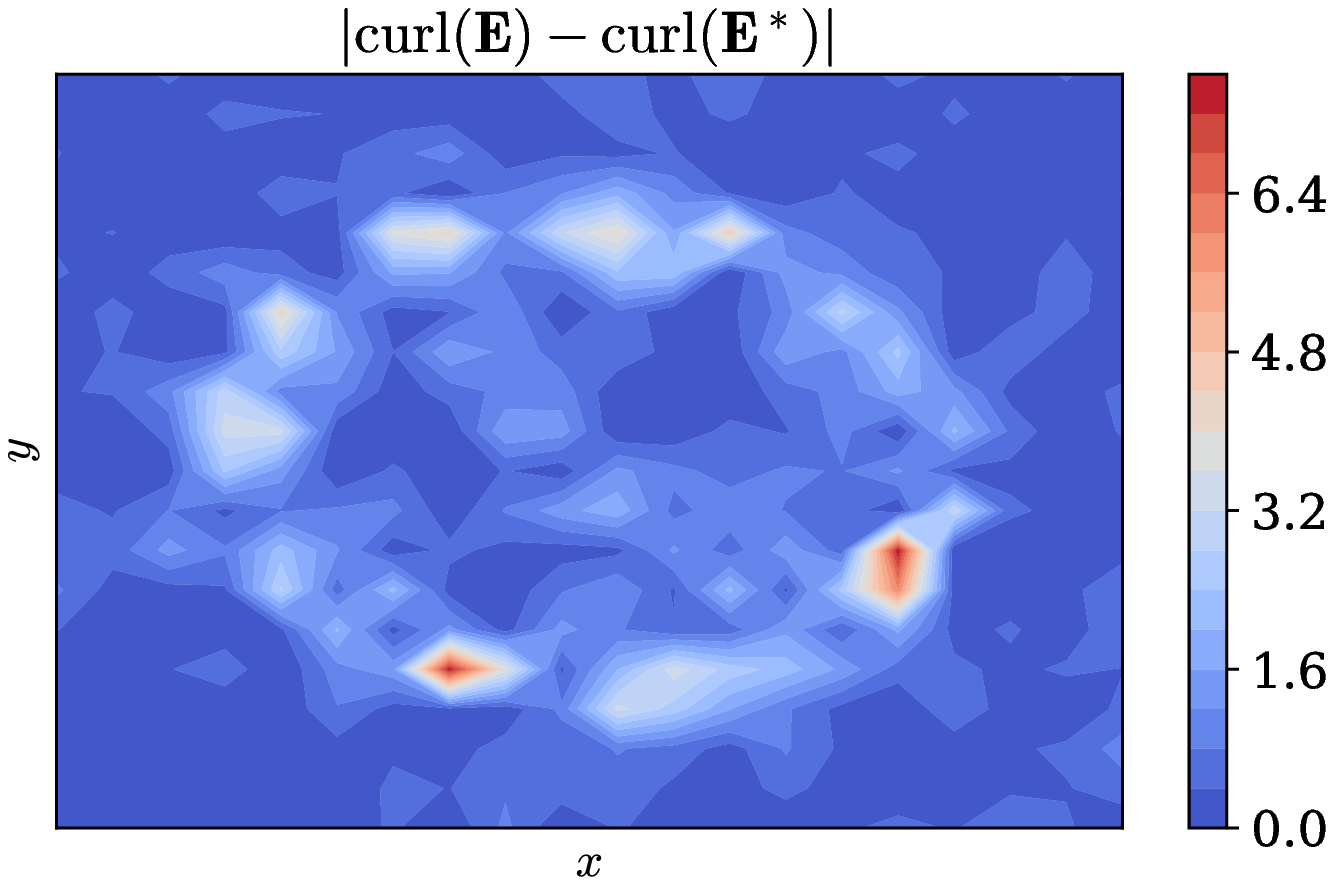}
		\caption{Error in the curl of the solution.}
		\label{fig:errorcurl_ex2}
	\end{subfigure}
	\caption{Solution and errors for the model Case~\hyperref[sec:case21]{2.1}.}
	\label{fig:P2_1}
\end{figure}

Figure~\ref{fig:P2_1} shows the obtained solution, the curl of the solution and the corresponding errors. As expected, the maximum errors are located close to the discontinuities of the parameters, which is inevitable in our implementation as we are using smooth NNs to approximate discontinuous functions. 

\begin{figure}[htbp!]
	\centering
	\includegraphics[width=0.55\textwidth]{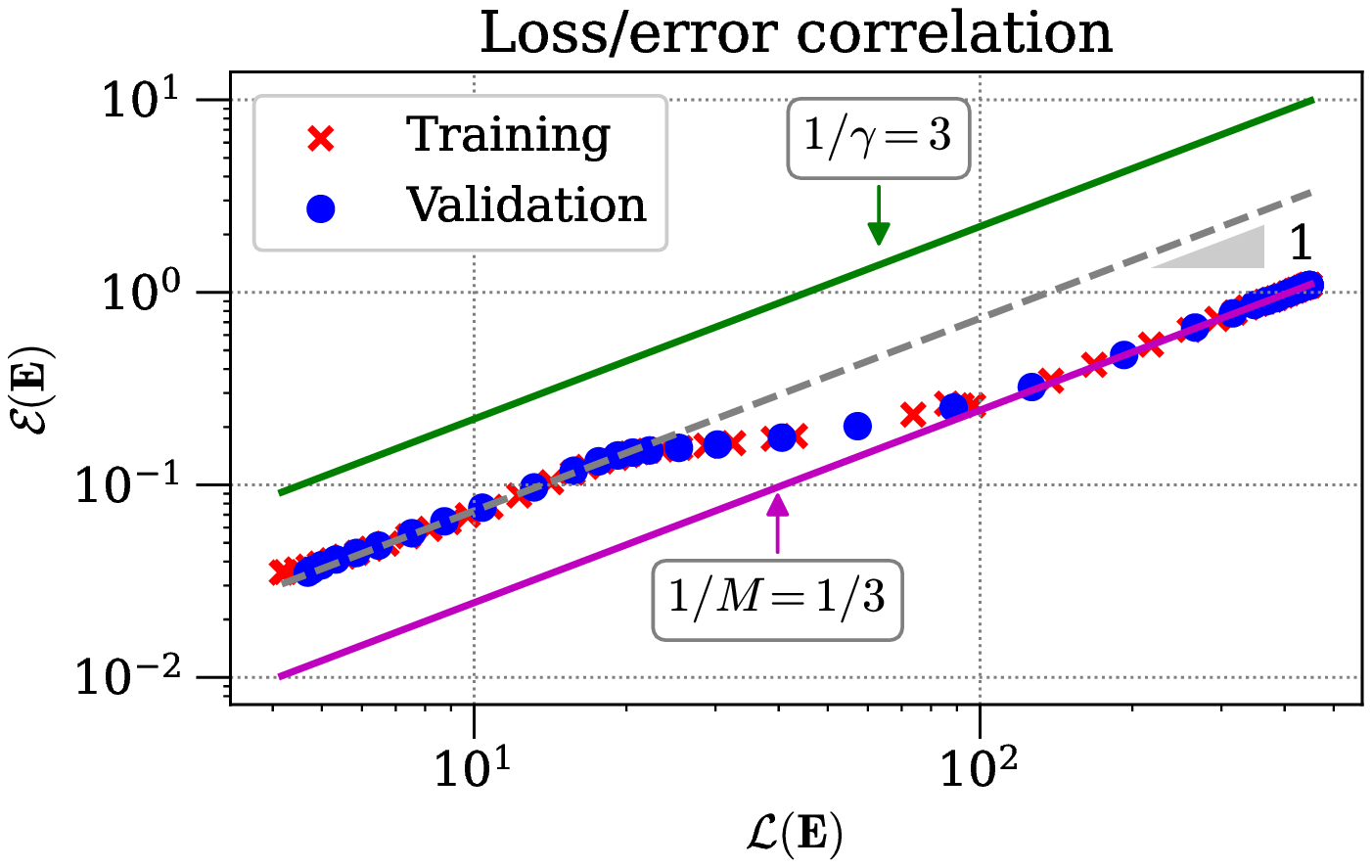}
	\vspace{-0.6cm}
	\caption{\centering The correlation between the loss $\LL(\tE)$ and the relative error in the solution $\mathcal{E}(\tE)$ during training and validation in Case~\hyperref[sec:case21]{2.1}.}
	\label{fig:P2_3}
\end{figure}

The relationship between losses and the relative error in the solution at each iteration is shown in Figure~\ref{fig:P2_3}. In the asymptotic regime, we obtain a linear relationship between the loss and the error. Since \eqref{eq:case2} relates with the inner product of $\Hc$, in this case we can use the Riesz representation theorem to estimate the equivalence constants $M$ and $\gamma$ in \eqref{eq:bounds}. A simple calculation give us $\frac{1}{M}=\frac{1}{3}$ and $\frac{1}{\gamma} =3$. Figure~\ref{fig:P2_3} exhibits the expected oscillation between these two parallel lines.

\begin{figure}[htbp!]
	\begin{subfigure}[b]{0.48\textwidth}
		\centering
		\includegraphics[width=1\textwidth]{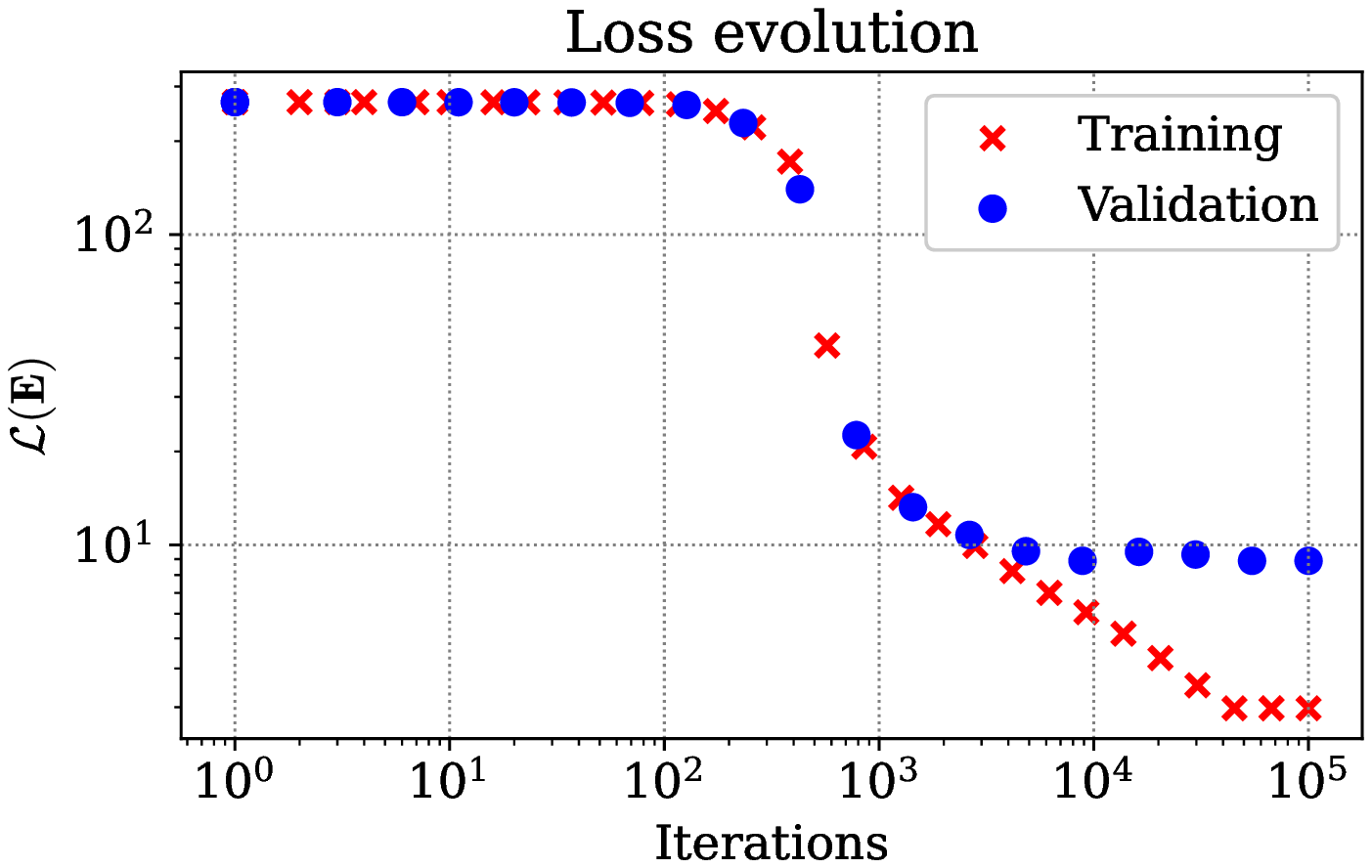}
		\caption{ \centering The evolution of the loss $\LL(\tE)$ on both the training and validation data sets.}
		\label{fig:lossevbada}
	\end{subfigure}
	\begin{subfigure}[b]{0.48\textwidth}
		\centering
		\includegraphics[width=1\textwidth]{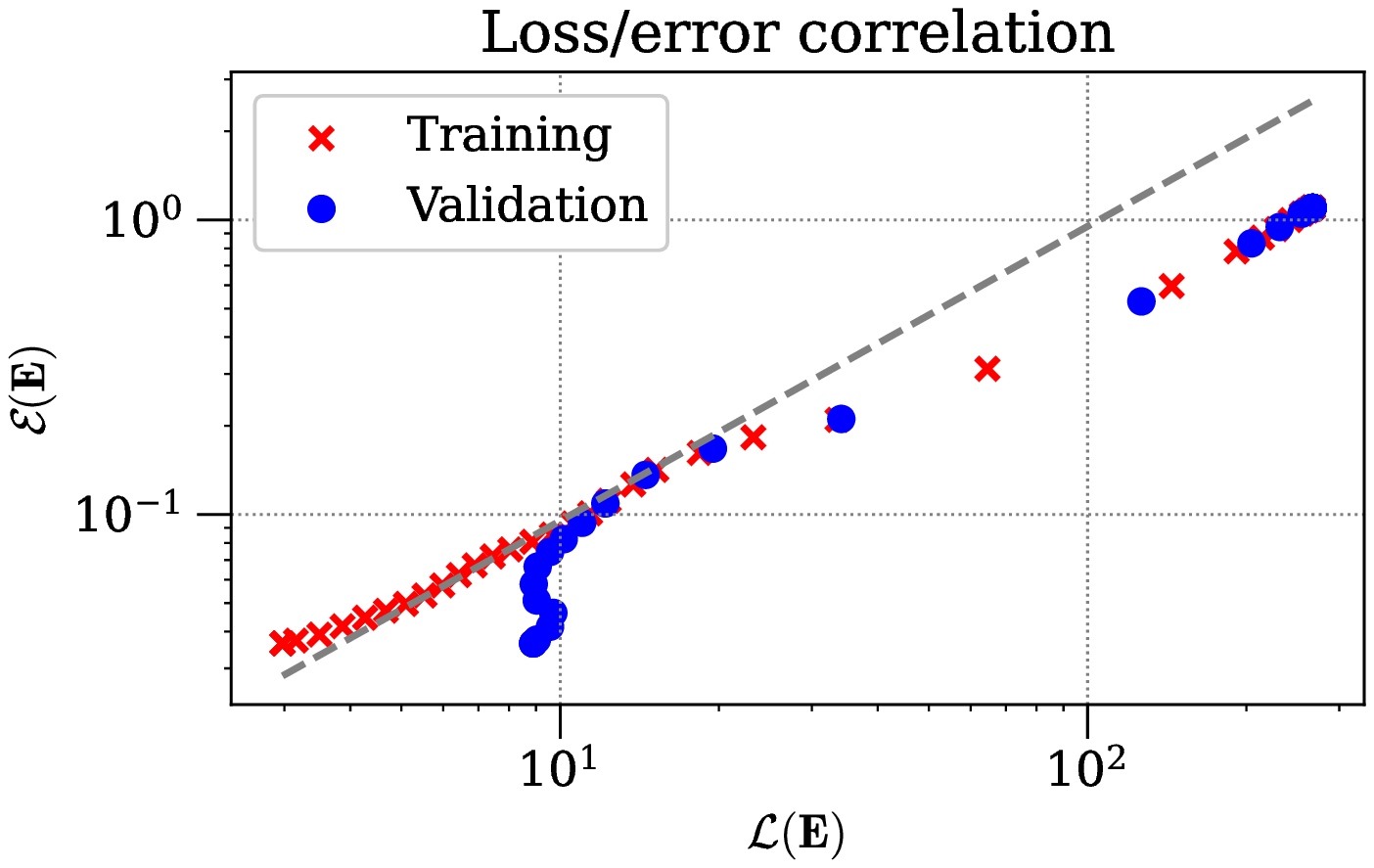}
		\caption{\centering The correlation between the loss and the relative error.}
		\label{fig:lossevbadb}
	\end{subfigure}
	\caption{\centering The evolution of the loss  and the correlation between the loss $\LL(\tE)$ and the relative error in the solution $\mathcal{E}(\tE)$ during training and validation in Case~\hyperref[sec:case21]{2.1} when using an inaccurate integration rule.}
	\label{fig:sol2_bad}
\end{figure}

Now, we remark the importance of accurately calculating the integrals in \eqref{eq:discloss}. Figure~\ref{fig:sol2_bad} shows the loss evolution and the correlation between the $H$(curl)-norm of the error and the loss when the number of integration points is equal to the number of integration points in the Case~\hyperref[sec:case21]{2.1}. There we use $N=M=100$ integration points for training the NN and $120$ points on each direction for validation. In Figure~\ref{fig:lossevbada} we observe a divergence in the training and validation losses after roughly $2000$ iterations, owing to integration errors, highlighting the need for accurate integration when solutions are of lower regularity. 

\subsection{Case 2.2. The physical variational form}\label{sec:case22}
Now, we modify  Case~\hyperref[sec:case20]{2} and seek for $\tE \in \Hcz$ satisfying
\begin{equation*}
	\int_\Om \mu^{-1}  \crl(\tE) \cdot \crl(\bphi)- \omega^2 \epsilon \tE \cdot \bphi \, d\tx = \int_\Om \tilde{\tJ} \cdot \bphi \, d\tx,
\end{equation*}
with $\omega = 1.25$. This modification corresponds to the weak formulation of the time-harmonic Maxwell's equations in \eqref{eq:weakMax} with positive discontinuous parameters $\mu$ and $\epsilon$.

\begin{figure}[htbp!]
	\begin{subfigure}[b]{0.48\textwidth}
		\centering
		\includegraphics[width=1\textwidth]{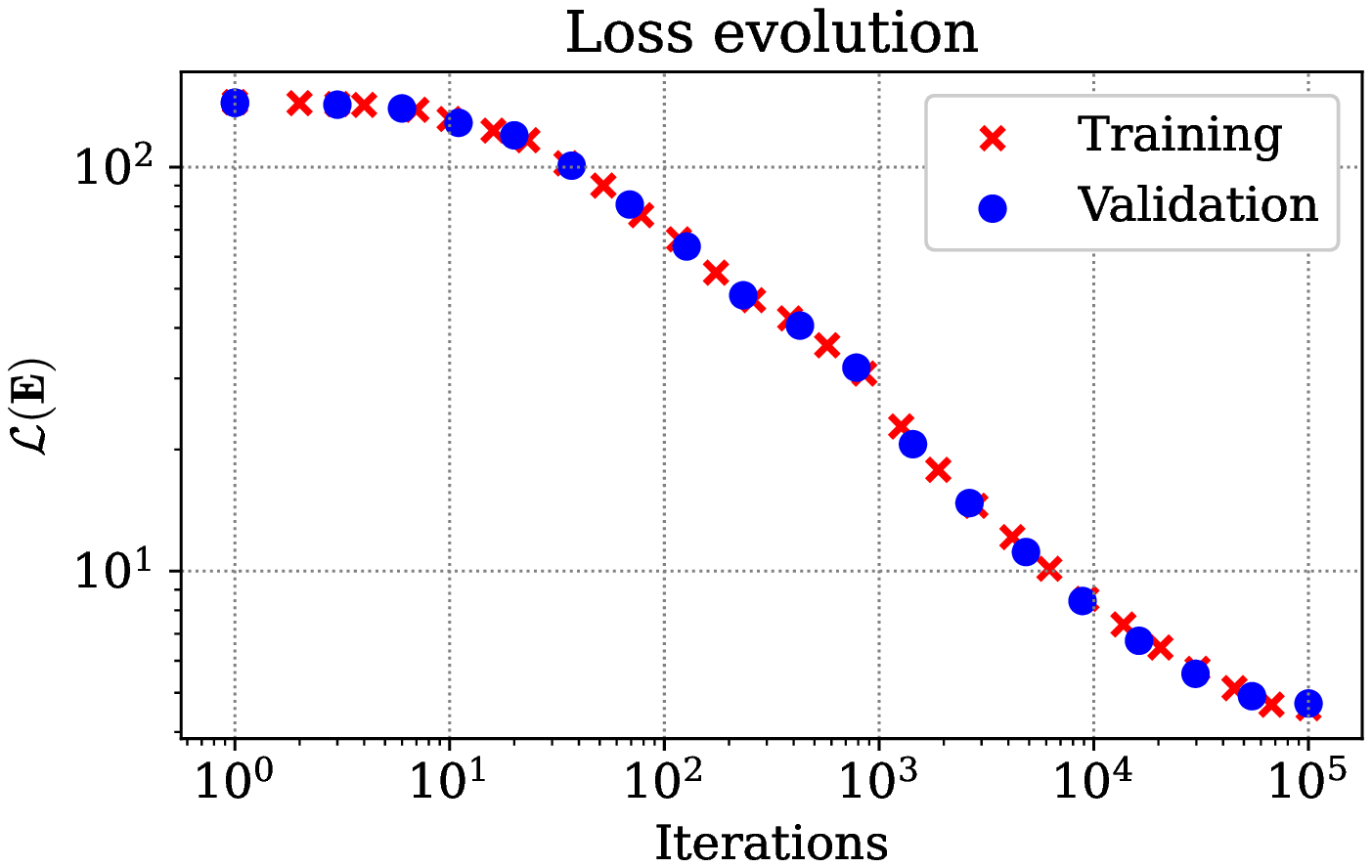}
		\caption{ \centering The evolution of the loss $\LL(\tE)$ on both the training and validation data sets. }
		\label{fig:lossev2}
	\end{subfigure}
	\begin{subfigure}[b]{0.48\textwidth}
		\centering
		\includegraphics[width=1\textwidth]{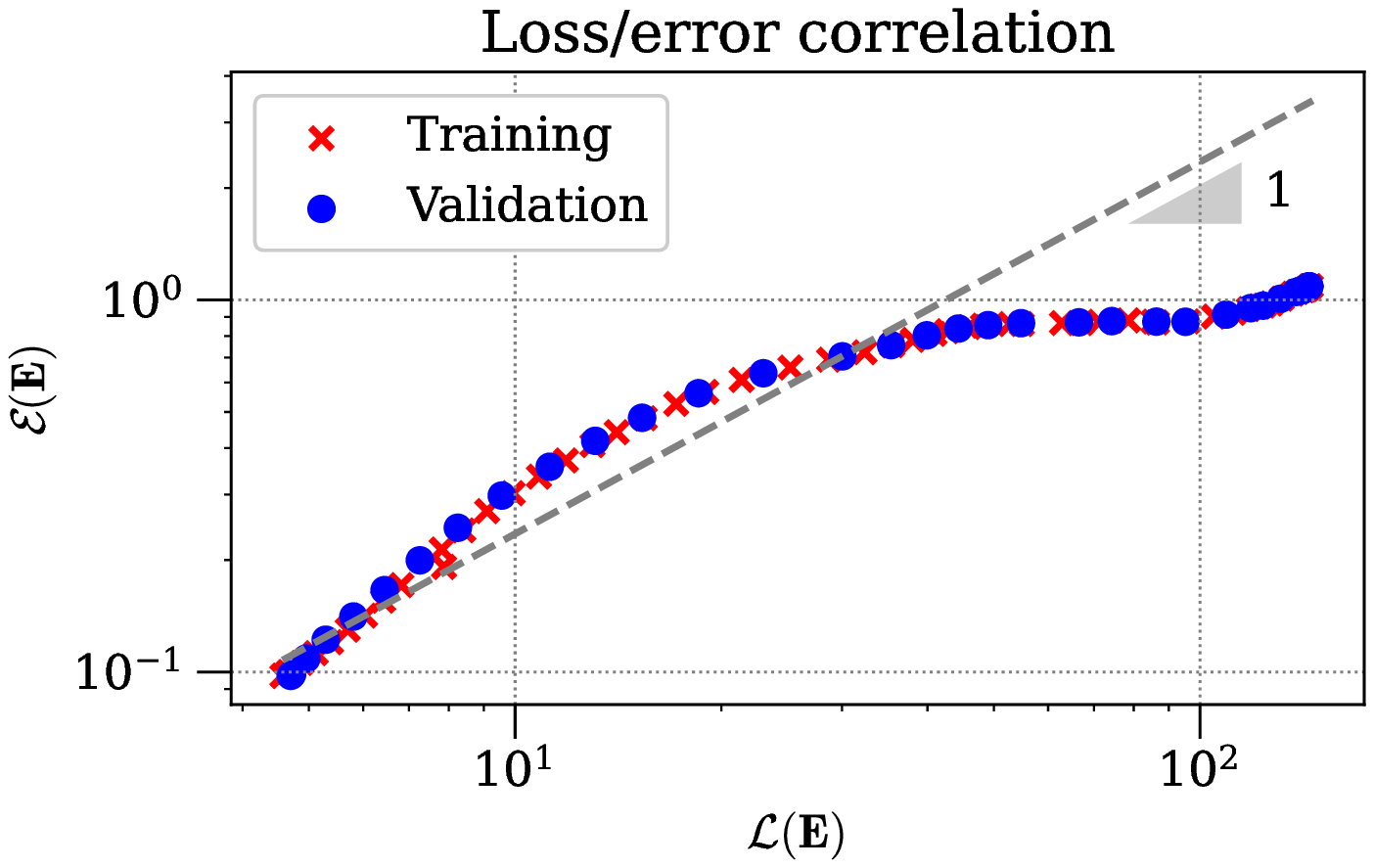}
		\caption{\centering The correlation between the loss and the relative error.}
		\label{fig:lossev}
	\end{subfigure}
	\caption{\centering The evolution of the loss and the correlation between the loss $\LL(\tE)$ and the relative error in the solution $\mathcal{E}(\tE)$ during training and validation in Case~\hyperref[sec:case22]{2.2}.}
	\label{fig:case2.2}
\end{figure}

In contrast to the Case~\hyperref[sec:case21]{2.1}, in Figure~\ref{fig:case2.2} one notices that the evolution of the loss decreases slower and the relation between the loss and the error is not linear. This effect is a clear consequence of the inclusion of the frequency $\omega$ that modifies the bounds of the error in \eqref{eq:bounds}.

\subsection{Case 3. Smooth solution in 3D}\label{sec:case3}

Take $\Om =[0,\pi]^3$. We consider the variational form: find $\tE \in \Hcz$ satisfying
\begin{equation}\label{eq:case3}
	\int_\Om \mu^{-1}  \crl(\tE) \cdot \crl(\bphi)- \omega^2 \epsilon \tE \cdot \bphi \, d\tx = \int_\Om \tilde{\tJ} \cdot \bphi \, d\tx \qquad \forall  \bphi \in \Hcz.
\end{equation}

We choose $\mu$ and $\epsilon$ to be constant equal to $1$ and $\omega= 1.5$.  Here, $\tilde{\tJ}$ is such that the exact solution is 
\begin{equation}
	\tE^*(x,y,z)= \begin{bmatrix} \sin(y)\sin(z)\sin(\omega x) \\ \sin(x)\sin(z)\sin(\omega y) \\ \sin(x)\sin(y)\sin(\omega z) \end{bmatrix}.
\end{equation}
 Here, we impose homogeneous-Dirichlet boundary conditions on $\partial\Omega$. We take a partition of $50$ integration points on each direction for training the NN and $60$ points for validation. Moreover, we use $50$ modes in both training and validation. 

\begin{figure}[htbp!]
	\centering
	\begin{subfigure}[b]{0.495\textwidth}
		\centering
		\includegraphics[width=1\textwidth]{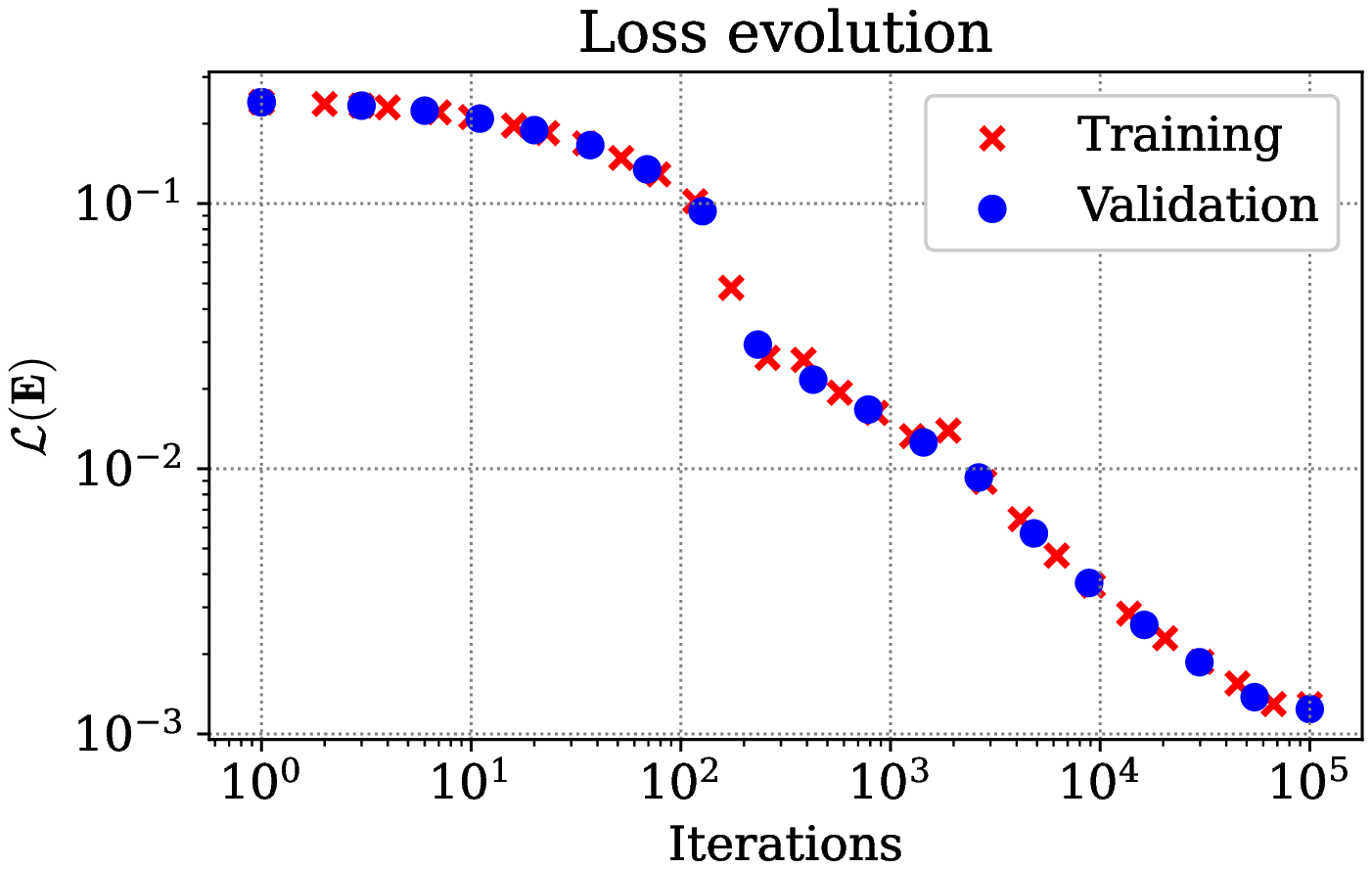}
		\caption{\centering The evolution of the loss $\LL(\tE)$ on both the training and validation data sets. }
		\label{fig:lossev_ex3}
	\end{subfigure}
	\begin{subfigure}[b]{0.495\textwidth}
		\centering
		\includegraphics[width=1\textwidth]{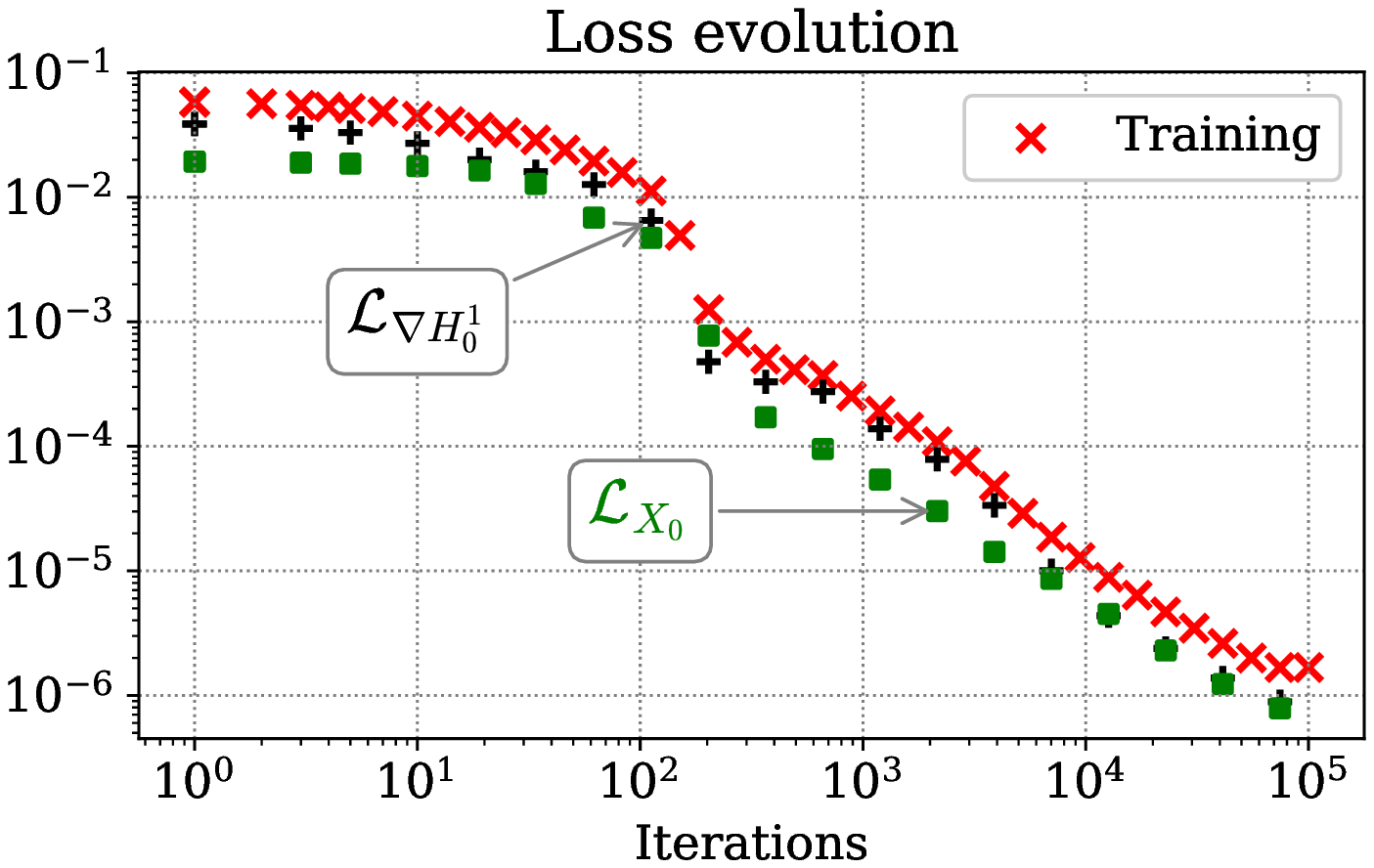}
		\caption{ \centering The contribution of the spaces $\nabla H_0^1(\Om)$ and $X_0(\Om)$ to the loss $\LL(\tE)^2$ on the training set.}
		\label{fig:lossev2_ex3}
	\end{subfigure}
	\caption{The evolution of the loss in Case~\hyperref[sec:case3]{3}.}
	\label{fig:P3_0}
\end{figure}

Similarly to the 2D cases, Figure~\ref{fig:lossev_ex3} shows the evolution of the loss on the training and validation data sets. After $10^{5}$ iterations, the losses stabilize and reach values upto $10^{-3}$. Figure~\ref{fig:lossev2_ex3} shows the contribution of the spaces $\nabla H_0^1(\Om)$ and $X_0(\Om)$ to the squared loss on the training set. 

\begin{figure}[htbp!]
	\centering
	\includegraphics[width=0.55\textwidth]{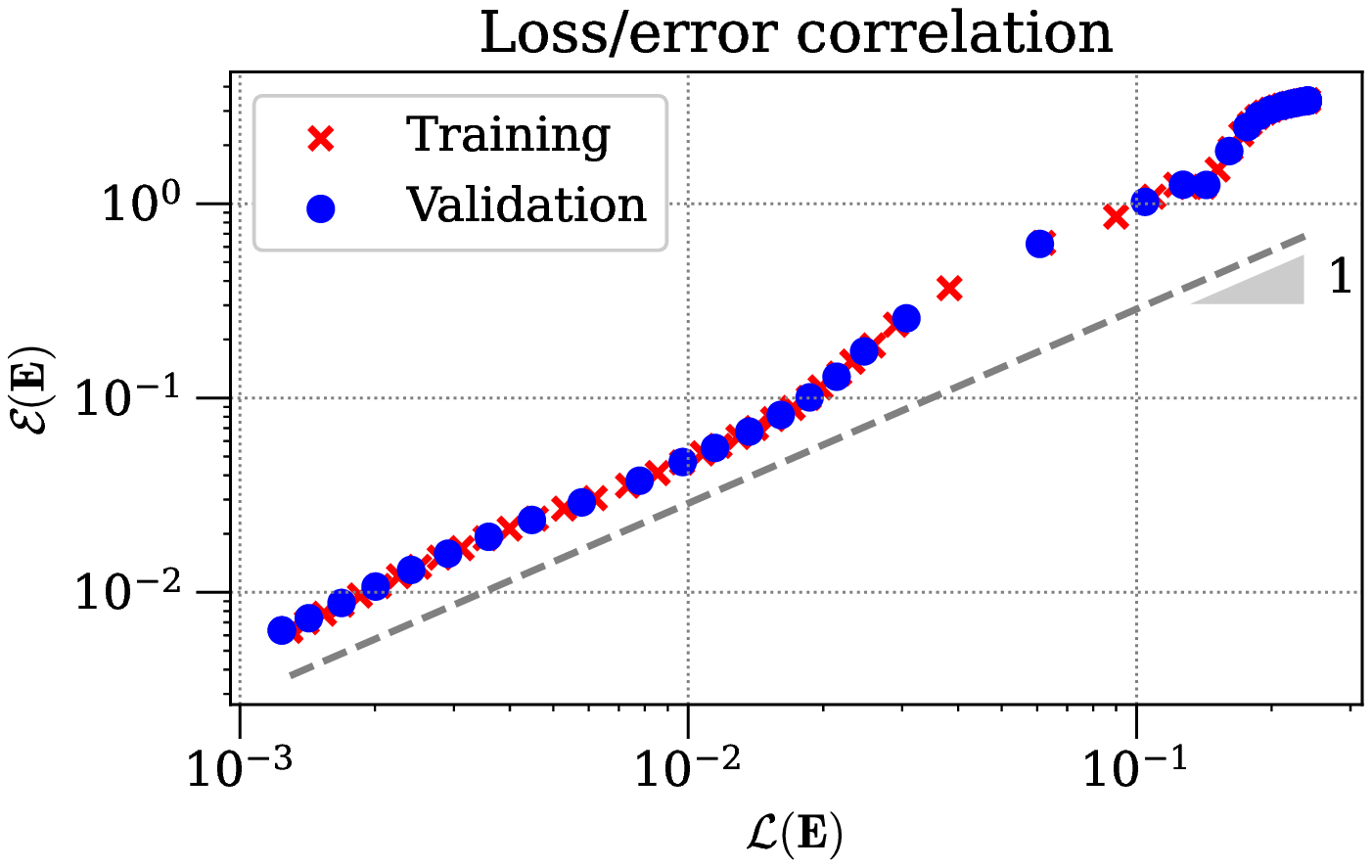}
	\vspace{-0.6cm}
	\caption{\centering The correlation between the loss $\LL(\tE)$ and the relative error in the solution $\mathcal{E}(\tE)$ during training and validation in Case~\hyperref[sec:case3]{3}.}
	\label{fig:P3_1}
\end{figure}

Finally, we illustrate the correlation between the loss and the relative error in the solution during training and validation in Figure~\ref{fig:P3_1}. In the asymptotic regime, the relationship between the loss and the error is linear.

\section{Conclusions}\label{sec:concl}

We extended the principles in  \cite{taylor2023deep} and implemented the DFR method for solving Maxwell's equations. In this case, we select the dual norm of the residual as the loss function of a NN solving Maxwell's problem. We rely on the weak formulation of Maxwell's electric field and the Helmholtz decomposition of the space $\Hcz$. We proposed orthonormal basis functions for each sub-space and used them to construct a computable loss function. To lower the computational cost of estimating integrals, we apply DST/DCT transformations in our discretized loss function. 
In two and three spatial dimensions, the numerical examples show a linear association between the loss and the $H$(curl)-norm of the error. We included examples with discontinuous parameters an noticed the importance of avoiding overfitting and integration errors.

We note that the DFR discussed here suffers the curse of dimensionality and that the implementation differs in 2D and 3D. Additionally, the DFR in 2D can only be applied to general domains if the eigenbasis of the Laplacian on those geometries is known. The 3D version is restricted to Cartesian product domains, which are rectangular or cylindrical domains based on the Laplacian's eigenbasis in 2D, which again must be known in order to implement the method. Future work will look at the method's scalability and investigate the use of subdomain-based local test functions. By doing so, we will consider more general geometries.

\section{Acknowledgements}

Jamie M. Taylor is supported by the Basque Government through the BERC 2018-2021 program and by the Spanish State Agency of Research through ``BCAM Severo Ochoa'' accreditation of excellence SEV-2017-0718 and through the project (PID2020-114189RB-I00 / AEI / 10.13039 / 501100011033). David Pardo has received funding from: the Spanish Ministry of Science and Innovation projects with references TED2021-132783B-I00, PID2019-108111RB-I00 (FEDER/AEI) and PDC2021-121093-I00 (MCIN / AEI / 10.13039/501100011033/Next Generation EU), the ``BCAM Severo Ochoa'' accreditation of excellence CEX2021-001142-S / MICIN / AEI / 10.13039/ 501100011033; the Spanish Ministry of Economic and Digital Transformation with Misiones Project IA4TES (MIA.2021.M04.008 / NextGenerationEU PRTR); and the Basque Government through the BERC 2022-2025 program, the Elkartek project SIGZE (KK-2021/00095), and the Consolidated Research Group MATHMODE (IT1456-22) given by the Department of Education. Ignacio Muga is supported by the Chilean National Agency for Research \& Development through the Fondecyt Project \#1230091.

 \bibliographystyle{siam}
 \bibliography{biblio}

\appendix
\addcontentsline{toc}{section}{Appendices}
\renewcommand{\thesubsection}{\Alph{section}.\Roman{subsection}}

\section{Basis for the space $\Hcz$} \label{sec:appA}
 This appendix contains classical results available, for example, in \cite{davies1995spectral, Evans1998}. First, in the following proposition, we construct a set of auxiliary eigenvectors of the operator $-\Delta$ in $H^1(\Om)$. 

 \begin{proposition}\label{prop1}
 	Let $\Om$ be a bounded, simply connected, and Lipschitz domain in $\mathbb{R}^n$, with $n = 2$ or $3$ and boundary $\partial \Om$. Assume that $\partial \Om$ consists of two disjoint parts $\Gamma_D$ and $\Gamma_N$, such that $\partial \Om = \Gamma_D \cup \Gamma_N$. Then, there exists a countable set of eigenfunctions $(\phi_k)_{k\in\I}$, with $\I$ being a set of indices, and corresponding eigenvalues $\lambda_k>0$ such that the following results hold, for all $k\in\I$
 	\begin{equation}\label{eq:lap3}
 		\begin{aligned}
 			-\Delta \phi_k & = \lambda_k \phi_k &\text{ in } \Om, \\
 			\phi_k & = 0 &\text{ on } \Gamma_D, \\
 			\nabla \phi_k \cdot \tn & = 0 &\text{ on } \Gamma_N.
 		\end{aligned}
 	\end{equation}
 	Moreover, the sequence $(\phi_k)_{k\in\I}$ forms an orthonormal basis for $L^2(\Om)$ and an orthogonal basis for $H^1(\Om)$.
 \end{proposition}
The proof of Proposition~\ref{prop1} is a classical result available, for example, in \cite[Theorem D.6.7]{Evans1998} and \cite[Theorem 1.8]{davies1995spectral}. Here, we only prove Propositions~\ref{propNabla} and \ref{propX2}.

\subsection*{Proof of Proposition~\ref{propNabla}}
\begin{proof}
	First we note that, since the space $\nabla H_0^1(\Om)$ as a subspace of $\Hc$ is isometric to $H_0^1(\Om)$, then the gradients of any orthogonal basis of $H^1_0(\Om)$ define an orthogonal basis of $\nabla H^1_0(\Om)$. 
	
	Take $\partial \Om = \Gamma_D$ in Proposition~\ref{prop1}. Then, there exists a set of eigenfunctions $(\phi_k)_{k\in\I}$ satisfying \eqref{eq:lap3}, which are called homogeneous-Dirichlet eigenvectors of $-\Delta$ in $\Om$. For each $k\in\I$, the function $\phi_k \in H^1_0(\Om)$, so the sequence $(\phi_k)_{k\in\I}$ forms an orthonormal basis for $L^2(\Om)$ and an orthogonal basis for $H^1_0(\Om)$. Therefore, the sequence $(\nabla \phi_k)_{k\in\I}$ forms an orthogonal basis of $\nabla H^1_0(\Om)$ due to the isometry between the spaces. The result in Proposition~\ref{propNabla} follows from the normalization of the functions $\nabla \phi_k$.
\end{proof}

\subsection*{Proof of Proposition~\ref{propX2}}
\begin{proof}
	Take $\partial \Om = \Gamma_N$ in Proposition~\ref{prop1}. Then, there exists a set of eigenfunctions $(\phi_k)_{k\in\I}$ satisfying \eqref{eq:lap3}. These are called homogeneous-Neumann eigenvectors of $-\Delta$ in $\Om$. We disregard the constant functions and define $\hat{\bpsi}_k  := \crl^*(\phi_k)$ for all $k$ in the index set $\I$. Each function $\phi_k$ is smooth, so one certainly has that $\hat{\bpsi}_k \in \Hc$ for all $k\in\I$. 
	
	From \eqref{eq:lap3}, we obtain that $\crl^*(\phi_k)\times\tn = \nabla\phi_k\cdot \tn = 0$ on $\Gamma_N$, so that $\hat{\bpsi}_k\in \Hcz$. Moreover, for all $u \in H_0^1(\Om)$ we have
	\[ ( \hat{\bpsi}_k,\nabla u )_{\Hc} = \int_{\Om} \crl^*(\phi_k) \cdot \nabla u \, d\tx = \int_{\Om} \phi_k \, \crl(\nabla u) \, d\tx = 0, \]
	meaning $\dv(\hat{\bpsi}_k) = 0$ weakly, so $\hat{\bpsi}_k\in X_0(\Om)$ for all $k\in\I$.
	
	To show that the sequence $(\hat{\bpsi}_k)_{k\in\I}$ forms an orthogonal basis of $X_0(\Om)$, we turn first to the orthogonality of the functions. Clearly, $\crl(\hat{\bpsi}_k) = -\Delta\phi_k = \lambda_k\phi_k$. Therefore, if $k \neq k'$, one obtains
	\begin{equation*}
		( \hat{\bpsi}_k,\hat{\bpsi}_{k'})_{\Hc} = \lambda_{k'}(\lambda_k+1)\int_\Om \phi_{k} \phi_{k'} \, d\tx = 0,
	\end{equation*} 
	as $\phi_k$ and $\phi_{k'}$ are orthogonal in $L^2(\Om)$. 
	
	Now, to show that they form a complete basis, we assume otherwise for the sake of contradiction, that there exists some non-zero $\tv \in X_0(\Om)$ such that $( \hat{\bpsi}_k,\tv)_{\Hc} = 0$ for all $k\in\I$. In this case, we can write that 
	\begin{equation*}
		\begin{split}
			0 = &( \hat{\bpsi}_k,\tv)_{\Hc} = \int_\Om (-\Delta\phi_k + \phi_k)\crl(\tv)\,d\tx = \left(\lambda_k+1\right) \int_\Om \phi_k\crl(\tv)\,d\tx. 
		\end{split}
	\end{equation*}
	As $\lambda_k>0$ and $(\phi_k)_{k\in\I}$ is a complete orthogonal basis for $L^2(\Om)$, disregarding constant functions, this implies that $\crl(\tv)$ is constant.  We write $\crl(\tv) = c \in \RR$. Accordingly, $\crl^*(\crl(\tv))=0$, so simple calculations lead to  \[c^2|\Omega|=\int_\Omega (\crl(\tv))^2\,dx = \int_\Omega \tv\cdot \crl^*(\crl(\tv))\,dx =0.\] Thus, $\crl(\tv) = 0$ and $\dv(\tv) = 0$ and this implies that $\tv = 0$, contradicting our original assumption.
	
	 We conclude that $(\hat{\bpsi}_k)_{k\in\I}$ is a complete orthogonal basis of $X_0(\Om)$. The result in Proposition~\ref{propX2} follows from the normalization of the functions $\hat{\bpsi}_k$.
\end{proof}

As discussed in Section~\ref{sec:basisF2}, to find an orthogonal basis for the space $X_0(\Om)$ in 3D is not as straightforward as in the 2D case. Nonetheless, in \cite{costabel2019maxwell} the authors demonstrate that the TE and TM modes constitute a complete and orthogonal basis of $X_0(\Om)$. Here, we only show that the functions $\bpsi_k^{\mathrm{TE}} \in X_0(\Om)$ and $\bpsi_k^{\mathrm{TM}} \in X_0(\Om)$ for all the appropriate indices $k$. 

\subsubsection*{TM modes}

Take $\Gamma_D = \overline{\Om^*} \times \partial I$ and $\Gamma_N = \partial \Om^*\times \bar{I}$ in Proposition~\ref{prop1}. Then, there exists a set of eigenfunctions $(p_k)_{k\in\I}$ satisfying \eqref{eq:lap3}. Disregarding the zero function, for all $k$ in the index set $\I$, we define $\hat{\bpsi}_k^{\mathrm{TM}}  := \crl(p_k \te)$, where $\te$ is the unit vector in the $z$-direction. The functions $\hat{\bpsi}_k^{\mathrm{TM}}$ are the $\crl$ of vector fields with a single non-zero component in the distinguished direction and they are divergence free.

Using the vector identities of the curl, the operator $\crl\text{-}\crl$ reduces to the negative Laplacian when acting on the vector field $\hat{\bpsi}_k^{\mathrm{TM}}$, i.e., 
\begin{equation}
	\crl (\crl (\hat{\bpsi}_k^{\mathrm{TM}})) = \nabla(\dv (\hat{\bpsi}_k^{\mathrm{TM}})) - \Delta \hat{\bpsi}_k^{\mathrm{TM}} =  - \crl(\Delta p_k \te) = \lambda_k\bpsi_k^{\mathrm{TM}},
\end{equation}
then $\hat{\bpsi}_k^{\mathrm{TM}}$ are eigenvectors of the curl-curl operator with corresponding eigenvalues $\lambda_k>0$. Direct calculations verify that $\hat{\bpsi}_k^{\mathrm{TM}}$ satisfies the boundary condition $\hat{\bpsi}_k^{\mathrm{TM}} \times\tn = 0$, implying that $\hat{\bpsi}_k^{\mathrm{TM}} \in X_0(\Om)$ for all $k\in\I$. 

\subsubsection*{TE Modes}

Similarly to the TM modes, take $\Gamma_N = \overline{\Om^*} \times \partial I$ and $\Gamma_D = \partial \Om^*\times \bar{I}$ in Proposition~\ref{prop1}. Then, there exists a set of eigenfunctions $(\phi_k)_{k\in\I}$ satisfying \eqref{eq:lap3}. Disregarding the zero function, for all $k$ in the index set $\I$, we define $\hat{\bpsi}_k^{\mathrm{TE}} := \crl(\crl(\phi_{k}\te))$, where $\te$ is the unit vector in the $z$-direction. The functions $\hat{\bpsi}_k^{\mathrm{TE}}$ are curl of smooth functions and divergence free. Moreover, $\hat{\bpsi}_k^{\mathrm{TE}}$ are eigenvectors of the curl-curl operator, which is inherited from the fact that $\phi_k$ are eigenvectors of the Laplacian. A direct calculation verifies that these satisfy the boundary condition $\hat{\bpsi}_k^{\mathrm{TE}}\times\tn = 0$, implying that $\hat{\bpsi}_k^{\mathrm{TE}} \in X_0(\Om)$ for all $k\in\I$.

\section{Basis functions on $\Om = [0,\pi]^n$}\label{sec:appB}

Here, we consider the domain $\Om = [0,\pi]^n$ with $n = 2$ or $3$. Table~\ref{tab1} shows the eigenvectors of the Laplacian in two and three dimensions for different types of boundary conditions. In 3D, we take $\Om = [0,\pi]^3 = [0,\pi] \times \Om^*$, where $\Om^* = [0,\pi]^2$. 

We note that the eigenvectors for the Laplacian on rectangular domains with distinct side lengths are readily obtained by rescaling those in Table~\ref{tab1}. The constructions of basis functions for $\Hcz$ on rectangular domains can be obtained via the same constructions outlined in Sections~\ref{secB1}-\ref{secB4}. The appearance of numerous multiplicative constants makes the arithmetic unwieldy, albeit simple, thus for simplicity we only include the calculations for $n$-dimensional cubes. 

In Table~\ref{tab1}, the mixed boundary conditions of type 1 are homogeneous-Dirichlet on $\{0,\pi\}\times \overline{\Om^*}$ and homogeneous-Neumann on $[0,\pi] \times \partial \Om^* $, and, the type 2 are the reversed. We let $\I$ be the set of indices $k$ such that $k = (k_1,k_2)$ in 2D, $k = (k_1,k_2,k_3)$ in 3D, and each $\phi_k$ satisfies the problem \eqref{eq:lap3}. We disregard some indices $k$, as shown in Table~\ref{tab1}. Notice that the eigenvalues corresponding to all the eigenvectors $\phi_k$ in Table~\ref{tab1} are $\lambda_k = |k|^2$. More details about the properties of these functions can be found in \cite{taylor2023deep}. 

\begin{table}[h]
	\begin{center}
		\renewcommand{\arraystretch}{1.8}
		{\small \begin{tabular}{|M{0.06\textwidth}|p{0.13\textwidth}|M{0.35\textwidth}|m{0.4\textwidth}|}
				\hline
				& Boundary conditions & $\phi_k$ & Index $k$ \\
				\hline
				2D & Dirichlet & $\frac{2}{\pi}\sin(k_1x)\sin(k_2y)$ & $k_\ii\in \ZZ_{\geq 0}$  \\
				\hline
				2D & Neumann & $\frac{2}{\pi}\cos(k_1x)\cos(k_2y)$ & { $k_\ii\in \ZZ_{\geq 0}$ with $k_{1} = 0$ xor $k_{2} = 0$} \\
				\hline
				3D &Dirichlet & $\frac{2\sqrt{2}}{\pi^{3/2}}\sin(k_1x)\sin(k_2y)\sin(k_3z)$ & $k_\ii\in \ZZ_{\geq 0}$ \\
				\hline
				3D &Mixed type 1 & $\sin(k_1x)\cos(k_2y)\cos(k_3z)$ & { $k_\ii\in \ZZ_{\geq 0}$ with $k_1>0$ and $k_{2} = 0$ xor $k_{3} = 0$} \\
				\hline
				3D & Mixed type 2 & $\cos(k_1x)\sin(k_2y)\sin(k_3z)$ & { $k_\ii\in \ZZ_{\geq 0}$ with $k_{2}>0$ and $k_{3}>0$} \\
				\hline
		\end{tabular}}
	\end{center}
	\vspace{-12pt}
	\caption{The eigenvectors of $-\Delta$ on $n$-dimensional cubes (2D and 3D) with different boundary conditions.}
	\label{tab1}
\end{table}

\subsection{Basis for $\nabla H^1_0(\Om)$ in 2D}\label{secB1}

Let $\Om = [0,\pi]^2$. From Table~\ref{tab1}, we take $\phi_k(x,y) = \frac{2}{\pi}\sin(k_1x)\sin(k_2y)$ satisfying \eqref{eq:lap3} with $\Gamma_D = \partial \Om$ and $k = (k_1,k_2)$. Here, $k_1$, $k_2$ and $k_3$ are integers greater or equal to zero, and we let $\I$ be the set of all such indices $k$. Notice that $ \| \nabla\phi_k \| _{[L^2(\Om)]^2} = |k|$ for all $k \in \I$. Using Proposition~\ref{propNabla}, we conclude that the set of functions $(\bphi_k)_{k\in \I}$, where each function $\bphi_k$ is defined as
\begin{equation*}
	\bphi_k  := \frac{1}{|k|}\nabla\phi_k = \frac{2}{\pi|k|}\begin{pmatrix}
		k_1\cos(k_1x)\sin(k_2y)\\
		k_2\sin(k_1x)\cos(k_2y)
	\end{pmatrix},
\end{equation*}
forms an orthonormal basis for $\nabla H^1_0(\Om) \subset \Hcz$ in 2D.

\subsection{Basis for $X_0(\Om)$ in 2D}\label{secB2}

Let $\Om = [0,\pi]^2$. From Table~\ref{tab1}, we take $\phi_k(x,y) = \frac{2}{\pi}\cos(k_1x)\cos(k_2y)$ satisfying \eqref{eq:lap3} with $\Gamma_N = \partial \Om$ and $k = (k_1,k_2,k_3)$. Here, $k_\ii\in \ZZ_{\geq 0}$ with $\ii = 1$ or $2$, and we let $\I$ be the set of all such indices $k$. We disregard the case when $k_1 = k_2 = 0$. Let $\hat{\bpsi}_k$ be the functions obtained by applying the operator $\crl^*$ to $\phi_k$, i.e., $\hat{\bpsi}_k  := \crl^*(\phi_k)$. A direct calculations lead to 
\begin{equation*}
	\| \hat{\bpsi}_k \| _{H(\crl,\Om)}^2 = \int_\Om |k|^4|\phi_k|^2+\left(\crl (\crl^*(\phi_k))\right)\phi_k \, d\tx = |k|^4+|k|^2,
\end{equation*}
as $\phi_k$ have unit $L^2$-norm. For all ${k\in \I}$, we obtain 
\begin{equation*}
	\bpsi_k  := \frac{\hat{\bpsi}_k}{\sqrt{|k|^4+|k|^2}} = \frac{c_k}{\sqrt{|k|^4+|k|^2}}\begin{pmatrix} k_2 \cos(k_1x)\sin(k_2y)\\ -k_1 \sin(k_1x)\cos(k_2y) \end{pmatrix}, 
\end{equation*} 
where 
\begin{equation*}
	c_k = \left\{\begin{array}{c c}
			\frac{2}{\pi} & k_1>0 \text{ and } k_2>0,\\
			\frac{\sqrt{2}}{\pi} & k_1 = 0 \text{ xor } k_2 = 0.
		\end{array}\right.
\end{equation*}
Using Proposition~\ref{propX2}, we conclude that the set $(\bpsi_k)_{k\in\I}$ forms an orthonormal basis of the space $X_0(\Om)$ in 2D. 

\subsection{Basis for $\nabla H^1_0(\Om)$ in 3D}\label{secB3}

Let $\Om = [0,\pi]^3$. From Table~\ref{tab1}, we take $\phi_k(x,y,z) = \frac{2\sqrt{2}}{\pi^{3/2}}\sin(k_1x)\sin(k_2y)\sin(k_3z)$ satisfying \eqref{eq:lap3} with $\Gamma_D = \partial \Om$ and $k = (k_1,k_2,k_3)$. Here, $k_\ii$ are integers greater or equal to zero, and we let $\I$ be the set of all such indices $k$. Notice that $ \| \nabla\phi_k \| _{[L^2(\Om)]^2} = |k|$ for all $k \in \I$. Using Proposition~\ref{propNabla}, we conclude that the set of functions $(\bphi_k)_{k\in \I}$, where each function $\bphi_k$ is defined as
\begin{equation*}
	\bphi_k  := \frac{1}{|k|} \nabla\phi_k = \frac{2\sqrt{2}}{ \pi^{3/2}|k|}\begin{pmatrix}
		k_1\cos(k_1x)\sin(k_2y)\sin(k_3z)\\
		k_2\sin(k_1x)\cos(k_2y)\sin(k_3z)\\
		k_3\sin(k_1x)\sin(k_2y)\cos(k_3z)
	\end{pmatrix},
\end{equation*}
forms an orthonormal basis for $\nabla H^1_0(\Om)$ in 3D.

\subsection{Basis for $X_0(\Om)$ in 3D}\label{secB4}
In cubic domains, i.e., $\Om = [0,\pi]^3$, the construction of the TM and TE modes is independent of the direction. Here, we take the $x$ direction (with unit vector $\te_1$) to be the distinguished direction for the TM and TE modes.

\begin{itemize}
	\item From Table~\ref{tab1}, we take $\phi_{k}(x,y,z) = \sin(k_1x)\cos(k_2y)\cos(k_3z)$ satisfying \eqref{eq:lap3} with $\Gamma_D = \{0,\pi\}\times \overline{\Om^*}$, $\Gamma_N = [0,\pi] \times \partial \Om^*$ and $k = (k_1,k_2,k_3)$. Here, $k_\ii\in \ZZ_{\geq 0}$ with $\ii = 1,2$ or $3$, and we let $\I$ be the set of all such indices $k$. We disregard the cases where $k_1 = 0$ or $k_2 = k_3 = 0$. 
	
	For all $k\in\I$, we define $\hat{\bpsi}_k^{\mathrm{TM}} := \crl(\phi_k \te_1)$. So, a direct calculation shows that \[\| \hat{\bpsi}_k^{\mathrm{TM}}\|_{\Hc}^2 = c'_k(1+|k|^2)(k_2^2+k_3^2),\] where $c'_k = \frac{\pi^3}{8}$ if $k_2>0$ and $k_3>0$, and $c'_k = \frac{\pi^3}{4}$ if $k_2 = 0$ or $k_3 = 0$. Thus, for all $k\in \I$, the normalized TM modes are 
	\begin{equation}\label{eq:TMmode}
		\bpsi_k^{\mathrm{TM}}  := \frac{1}{\sqrt{c'_k(1+|k|^2)(k_2^2+k_3^2)}}\begin{pmatrix}
			0 \\ -k_3\sin(k_1x)\cos(k_2y)\sin(k_3z)\\ k_2\sin(k_1x)\sin(k_2y)\cos(k_3z)
		\end{pmatrix}, \forall k\in\I.
	\end{equation}

	\item From Table~\ref{tab1}, we take $\phi_{k}(x,y,z) = \cos(k_1x)\sin(k_2y)\sin(k_3z)$ satisfying \eqref{eq:lap3} with $\Gamma_N = \{0,\pi\}\times \overline{\Om^*}$, $\Gamma_D = [0,\pi] \times \partial \Om^*$ and $k = (k_1,k_2,k_3)$. Here, $k_\ii\in \ZZ_{\geq 0}$ with $\ii = 1,2$ or $3$, and we let $\I$ be the set of all such indices $k$. We disregard the cases where $k_2 = 0$ or $k_3 = 0$. 

For all $k\in \I$, we define $\hat{\bpsi}_{k}^\mathrm{TE} = \crl(\crl(\phi_{k}\te_1))$. The norm of $\hat{\bpsi}_{k}^\mathrm{TE}$ in $H(\crl,\Om)$ is \[\|\hat{\bpsi}_{k}^\mathrm{TE}\|_{H(\crl,\Om)}^2 = \frac{\pi^3}{8}|k|^2(1+|k|^2)\left(|k|^2-k_1^2\right).\] This gives the TE modes as 
\begin{equation}\label{eq:TEmode}
	\bpsi_k^{\mathrm{TE}} = \frac{2\sqrt{2}}{\pi^{3/2}\sqrt{(1+|k|^2)|k|^2\left(|k|^2-k_1^2\right)}}
	\begin{pmatrix}
		(k_2^2+k_3^2)\cos(k_1x)\sin(k_2y)\sin(k_3z)\\
		-k_1k_2\sin(k_1x)\cos(k_2y)\sin(k_3z)\\
		-k_1k_3\sin(k_1x)\sin(k_2y)\cos(k_3z)
	\end{pmatrix}, \forall k\in\I.
\end{equation}
\end{itemize}
	
 \end{document}